\documentclass{amsart}
\usepackage{tikz}
\usepackage{graphicx}
\usepackage{enumerate,url}
\newtheorem{theo}{Theorem}

\newtheorem{lemma}[theo]{Lemma}
\newtheorem{prop}[theo]{Proposition}
\newtheorem{cor}[theo]{Corollary}

\newtheorem{assums}[theo]{Assumptions}

\newtheorem{exa}[theo]{Example}

\newtheorem{rems}[theo]{Remarks}
\newtheorem{rem}[theo]{Remark}
\numberwithin{equation}{section}
\numberwithin{theo}{section}
\title[Vector-valued equations with coupled dynamic boundary conditions]{Vector-valued heat equations and networks\\ with coupled dynamic boundary conditions}
\keywords{Dynamic boundary conditions, invariance properties, coupled boundary conditions for vector-valued diffusion, strongly continuous semigroups}
\subjclass[2000]{47D06, 35K50}
\author{ Delio Mugnolo}
\address{Institut f\"ur Analysis, Universit\"at Ulm, Helmholtzstra{\ss}e 18, D-89081 Ulm, Germany}
\email{delio.mugnolo@uni-ulm.de}
\thanks{I warmly thank Robin Nittka (Ulm) for helpful suggestions.}

\begin{document}

\maketitle

\begin{abstract}
Motivated by diffusion processes on metric graphs and ramified spaces, we consider an abstract setting for interface problems with coupled dynamic boundary conditions belonging to a quite general class. Beside well-posedness, we discuss positivity, $L^\infty$-contractivity and further invariance properties. We show that the parabolic problem with dynamic boundary conditions enjoy these properties if and only if so does its counterpart with time-independent boundary conditions. Furthermore, we prove continuous dependence of the solution to the parabolic problem on the boundary conditions in the considered class. 
\end{abstract}

% 
% \bigskip
% 
% {\sc Resum\'e.}
% Motiv\'ees par des processus de diffusion sur des r{\'e}seaux m\'etriques et des espaces ramifi{\'e}s nous consid\'erons une approche abstraite aux problemes aux interfaces avec des conditions aux limites dynamiques et coupl\'ees. Nous \'etudions le caract\`ere bien pos\'e, la positivit\'e, la contractivit\'e dans $L^\infty$ et des propri\'etes qualitatives additionelles. Nous d\'emontrons que le probl{\`e}me parabolique avec des conditions aux limites dynamiques satisfait ces propri\'et\'es si et seulement si le probl\`eme parabolique avec des conditions aux limites \emph{non} dynamiques le satisfait. Nous prouvons que la d\'ependance de la solution du probl\`eme parabolique  aux conditions aux limites dans la classe consider\'ee est continue.

\section{Introduction}

Elliptic systems with coupled boundary conditions have been attracting broad attention at least since~\cite{AgmDouNir64}. 
A classical approach is based on interpreting interface conditions of an elliptic system as boundary conditions of a vector-valued elliptic equation. This leads to introducing differential operators acting on spaces of vector-valued functions. A parabolic theory for this kind of operators has been recently developed, see e.g.~\cite{Ama01,DenHiePru03}.

A particularly interesting application of the theory of elliptic systems is given by so-called \emph{networks} and \emph{quantum graphs}, see e.g.~\cite{AliBelNic01,Kuc08} and references therein. Their generalisation to $n$-dimensional problems has appeared already in~\cite{Lum80a}, where the related notion of \emph{ramified space} has been proposed. Having in mind applications to quantum graphs, Kuchment has proposed in~\cite{Kuc04} a class of coupled, time-independent boundary conditions for $1$-dimensional elliptic systems. Kuchment's formalism allows for a very efficient variational approach, but the tradeoff is that his boundary conditions are only a proper subset of those considered in~\cite{AgmDouNir64} -- or, in the specific context of quantum graphs, in~\cite{KosSch99}. However, it is remarkable that Kuchment's conditions give rise exactly to \emph{all self-adjoint realisations} of the Schr\"odinger operator on a metric graph, under a mild locality assumption.

In the companion paper~\cite{CarMug09}, Cardanobile and the author have generalized Kuchment's formalism to the case of $n$-dimensional vector-valued diffusion and characterized several properties of the parabolic problem in dependence on the chosen boundary conditions. The aim of this paper is to provide the extension of the theory in~\cite{CarMug09} to the case of \emph{dynamic} boundary conditions of Wentzell--Robin-type.

Although we are soon going to consider the general case, let us start by briefly focusing on the $1$-dimensional setting of \emph{networks} (or \emph{quantum graphs}).

\begin{exa}\label{firstexa}
Let $N\in\mathbb N$ and consider the prototypical case of a diffusion problem
\begin{equation}\tag{TDPS}
\left\{
\begin{array}{rcll}
\dot{u}_j(t,x)&=& u_j''(t,x), &t\geq 0,\; x\in(0,\infty), \; j=1,\dots,N,\\
u_j(t,0)&=&u_\ell (t,0)=:\psi(t), &t\geq 0,\; j,\ell=1,\ldots,N,\\
\dot{\psi}(t)&=&\sum_{j=1}^N u'_j(t,0)&t\geq 0,
\end{array}
\right.
\end{equation}
on a metric graph -- more precisely, on a semi-infinite star with $N$ edges $e_1,\ldots,e_N$ on whose center a dynamic Kirchhoff-type boundary condition is imposed along with a standard continuity assumption. Each edge is parametrized as a $(0,\infty)$-interval, where $0$ is identified as the center of the star. Therefore, the function $u_j$ describing the diffusion on the edge $e_j$ maps $[0,\infty)\times [0,\infty)$ to $\mathbb C$, while $\psi:[0,\infty)\to\mathbb C$ describes the time evolution of the common boundary value in the center. It is known that the associated initial value problem is well-posed, as discussed, e.g. in~\cite{AliNic93,Bel94,MugRom07}. 

Laplace operators with dynamic boundary conditions appear as limiting cases of approximation schemes considered in~\cite{KucZen03,ExnPos05}. The cable model of a dendritical tree proposed by Rall in~\cite{Ral59} also leads to analogous network diffusion problems, cf.~\cite{Cam80,Nic85}: a thorough biomathematical investigation of them has been performed in a series of four papers beginning with~\cite{MajEvaJac93}.\\
\begin{center}
\begin{tikzpicture}[style=thick]
\foreach \x in {18,90,162,214,270,306}{
\draw[fill] (\x:0cm) -- (\x:1.8cm);
\draw[fill] (\x:1.85cm) -- (\x:1.9cm);
\draw[fill] (\x:1.95cm) -- (\x:2cm);
\draw[fill] (\x:2.05cm) -- (\x:2.1cm);
\draw[fill] (0:0cm) circle (2pt) node[above left]{$0$};
}
\end{tikzpicture}\\
\emph{A semi-infinite star with $6$ edges.}\\[1em] 
\end{center}

We can rephrase $\rm(TDPS)$ by considering the orthogonal projection $P_Y$ of ${\mathbb C}^N$ onto the subspace $Y:=\langle {\mathbf 1}\rangle$ spanned in ${\mathbb C}^N$ by the vector 
$${\mathbf 1}:=(1,\ldots,1).$$ 
Observe that the unknown can be thought of as a function $u:(0,\infty)\to{\mathbb C}^N$, so that the network diffusion problem simply becomes
$$\dot{u}(t,x)=u''(t,x),\qquad t\ge 0,\; x\in (0,\infty),$$
with suitable boundary conditions in $0$. More precisely, the continuity condition in the star's center -- given by the second equation in $\rm(TDPS)$ -- amounts to require that $u(t,0)\in \langle {\mathbf 1}\rangle$ for all $t\ge 0$, i.e.,
$$P_Y(u(t,0))=u(t,0),\qquad t\ge 0,$$
while the dynamic boundary condition equivalently reads 
$$\dot{u}(t,0)=P_Y(\dot{u}(t,0))=N P_Y(u'(t,0))=-P_Y\left(\frac{\partial u}{\partial \nu}(t,0)\right),\qquad t\ge 0.$$
Hence, the dynamic boundary condition is an equation living in the ($1$-dimensional) \emph{boundary space} $Y=\langle {\mathbf 1}\rangle$.

This kind of boundary conditions also arises in the mathematical modelling of string networks with masses at the nodes. They play an important role in the control theory of wave and beam equations: investigations in this direction go back at least to~\cite[{\S}2.7]{LagLeuSch94} and~\cite{HanZua95}.
\end{exa}

The goal of the present article is to generalize the setting discussed in the above example. Let $\Omega$ be a smooth open domain in ${\mathbb R}^n$ with boundary $\Gamma:=\partial\Omega$. Let $H$ be a separable complex Hilbert space. In particular, Bochner spaces $L^2(\Omega;H)$ and $L^2(\Gamma;H)$ become separable complex Hilbert spaces when endowed with the canonical scalar products
$$(f|g)_{L^2(\Omega;H)}:=\int_\Omega (f(x)|g(x))_H dx,\qquad f,g\in L^2(\Omega;H),$$
and
$$(f|g)_{L^2(\Gamma;H)}:=\int_\Gamma (f(z)|g(z))_H d\sigma(z),\qquad f,g\in L^2(\Gamma;H).$$
Let $\mathcal Y$ be a closed subspace of $L^2(\Gamma;H)$ and hence a Hilbert space in its own right with respect to the scalar product induced by $L^2(\Gamma;H)$. Vector-valued Sobolev spaces can be introduced recursively just like in the 
scalar-valued case. I.e., one first lets $H^0(\Omega;H):=L^2(\Omega;H)$, hence defines
\begin{equation}\label{sob1}
\begin{array}{ll}
&H^k(\Omega;H):=\left\{f\in H^{k-1}(\Omega;H):\exists \nabla f:=g\in L^2(\Omega;H^n)\hbox{ s.t.}\right.\\[2pt]
&\qquad\qquad\qquad\qquad\left. \int_\Omega f(x)\nabla h(x) dx=-\int_\Omega g(x) h(x) dx\hbox{ for all }h\in C^\infty_c(\Omega;{\mathbb C})\right\},\qquad k=1,2,\ldots,
\end{array}
\end{equation}
and finally introduces spaces of fractional order by standard complex interpolation. (Here we denote by $H^n$ the Hilbert space defined as the Cartesian product of $n$ copies of $H$.)
In particular, $H^1(\Omega;H)$ is a Hilbert space with respect to the scalar product
\begin{equation}\label{sob2}
(f|g)_{H^1(\Omega;H)}:=\int_\Omega (\nabla f(x)|\nabla g(x))_{H^n} dx + \int_\Omega (f(x)|g(x))_{H} dx,\qquad f,g\in H^1(\Omega;H).
\end{equation}
We emphasize that vector-valued Sobolev spaces are introduced using scalar-valued test functions, hence the integral appearing in~\eqref{sob1} is vector-valued (i.e., a Bochner integral) whereas those appearing in~\eqref{sob2} are scalar-valued (i.e., Lebesgue integrals). It is well-known that the usual trace and normal derivative operators 
$$C(\overline{\Omega})\ni u\mapsto u_{|\Gamma}\in C(\Gamma)\qquad\hbox{and}\qquad C^1(\overline{\Omega})\ni u\mapsto \frac{\partial u}{\partial \nu}\in C(\Gamma)$$
extend to operators acting between Sobolev spaces of scalar-valued functions. In fact, they can be canonically defined in the vector-valued case, too -- e.g. by means of~\cite[Thm.~4.5.1]{Gra08}. With an abuse of notation we therefore denote by $u_{|\Gamma}$ and $\frac{\partial u}{\partial \nu}$ the trace and normal derivative (in the sense of distributions) of a function $u:\Omega\to H$.

We are now in the position to generalise the one-dimensional setting presented in Example~\ref{firstexa} by allowing for more general coupling conditions at the interface and consider the abstract boundary-value problem
\begin{equation}\tag{AS}
\left\{
\begin{array}{rcll}
\frac{\partial}{\partial t}{u}(t)&=& \Delta u(t), &t\geq 0,\\
u(t)_{|\Gamma}&\in& {\mathcal Y}, &t\geq 0,\\
\frac{\partial }{\partial t}u(t)_{|\Gamma}&=&-P_{\mathcal Y}\frac{\partial u(t)}{\partial\nu}, &t\geq 0,\\
\end{array}
\right.
\end{equation}
where $P_{\mathcal Y}$ denotes the orthogonal projection of $L^2(\Gamma;H)$ onto the closed subspace $\mathcal Y$.
In the $1$-dimensional case of finite quantum graphs, the investigation of such a problem has been sketched in~\cite[{\S}4]{KanKlaVoi09}.

\begin{exa}\label{dyenkirch}
Let $\Omega=(0,\infty)$ and $H={\mathbb C}^N$, so that $L^2(\Gamma)\equiv L^2(\{0\})\equiv {\mathbb C}^N$. Take 
$$\mathcal Y:=\langle({\mathbf 1})\rangle=\{c\in {\mathbb C}^N: c_1=\ldots=c_N\}.$$ 
Then
$$P_{\mathcal Y}=\frac{1}{N}\begin{pmatrix}
1 & \ldots & 1\\
\vdots &\ddots & \vdots\\
1 & \ldots & 1
   \end{pmatrix}
%\quad\hbox{and}\quad 
%P_{\mathcal Y^\perp}=\frac{1}{N}\begin{pmatrix}
%N-1 & & -1\\
% &\ddots & \\
%-1 & & N-1
%   \end{pmatrix}
$$
and one sees that $\rm(AS)$ is just a reformulation of $\rm(TDPS)$ considered in Example~\ref{firstexa}.
\end{exa}

\begin{exa}
Let again $\Omega=(0,\infty)$ and $H={\mathbb C}^N$. If $N=1$ and $\mathcal Y=L^2(\Gamma;H)\equiv{\mathbb C}$, then the first boundary condition in $\rm(AS)$ is void and $\rm(AS)$ is the reformulation of a scalar-valued heat equation with Wentzell--Robin boundary conditions, see e.g.\ the recent contributions in~\cite{FavGolGol02,AreMetPal03,MugRom06,VazVit08}. If instead $\mathcal Y=\{0\}$, then $\rm(AS)$ reduces to a heat equation with Dirichlet boundary conditions. For $N=1$, these are the only possible choices for $\mathcal Y$, but for $N\ge 2$ we have infinitely many new boundary conditions that in some sense interpolate between Dirichlet and Wentzell--Robin ones. This is crucial when setting up a Courant--Fischer min-max formula, cf.~\cite{BelMug10}. 
\end{exa}

\begin{exa}
For $H={\mathbb C}^N$ the elliptic problem with \textit{dynamic interface conditions} -- a vector-valued version of Wentzell--Robin boundary conditions -- has been considered in~\cite[{\S}III.4.5]{Sho94}. In~\cite{Nic88}, even more general elliptic interface problems have been considerd under the very general assumption that the given system can even consist of several metric spaces with different Hausdorff dimensions, see also~\cite{BelNic96}, see also~\cite{AliNic93}.

As already mentioned, the general case of a diffusion equation equipped with coupled (either dynamic or time-independent) boundary conditions is mostly motivated by the theories of quantum graphs and parabolic network equations, but it also appears in higher dimensional applications, in particular in biomathematical models -- see e.g.~\cite{MorFisPes08} and references therein. 
\end{exa}

In this article we restrict to the case of dynamic boundary conditions only. However, the general case of mixed dynamic/time-independent boundary conditions (typically appearing in models from the applied sciences, see e.g.~\cite{MajEvaJac93}) can be easily treated combining the results presented here and those from~\cite{CarMug09}. 

\bigskip
In Section~\ref{prelim} we introduce our abstract framework and deduce a well-posedness result. The above examples suggest that the vector-valued setting -- although equivalent to the that based on a network (or ramified space) formalism -- is more efficient. In fact, its flexibility allows to simply introduce whole families of spaces $\mathcal Y$. Consequently, completely new questions arise. For example, one may wonder how the solution to the heat equation with boundary conditions as in $\rm(AS)$ depends on $\mathcal Y$: it will be shown in Theorem~\ref{olafrev} that this dependence is continuous in norm under very natural assumptions. This result is interesting in that it does not have a scalar-valued \emph{pendant}. We also extend to the vector-valued case a result on continuous dependence on parameters obtained in the scalar-valued case in~\cite{CocFavGol08}. 

We consider invariance of order intervals and subspaces in Section~\ref{posdom}, showing in particular a tight relation beween the properties of the heat semigroup governing the problem with time-independent (i.e., Robin-type vector-valued) boundary conditions and its dynamic counterpart. We will observe some unexpected phenomena: e.g., the semigroups governing these diffusion problems are in general not submarkovian -- not even positivity preserving. 

To discuss these behaviours in detail, in Section~\ref{examplesect} we focus on the setting of Example~\ref{firstexa}. It turns out that even in the simple context of diffusion on a semi-infinite star with finitely many edges, unexpected dynamics arises after chosing appropriate boundary conditions.
%This leads, e.g., to the full characterisation of all local boundary conditions leading to a symmetric submarkovian heat semigroup on a dihedron in Example~\ref{diheall}. 

Finally, in Section~\ref{dynneum} we briefly discuss the general properties of a similar but different kind of dynamic boundary condition, where the normal derivative -- rather than the trace -- undergoes a time evolution.

\section{Preliminary results}\label{prelim}

To begin with, we make our standing assumptions precise. 

As in the previous section, let $H$ be a separable complex Hilbert space, $\Omega$ be an open domain in ${\mathbb R}^n$ with $C^1$ boundary $\Gamma:=\partial\Omega$ and $\mathcal Y$ be a closed subspace of $L^2(\Gamma;H)$. In the rest of the paper we are going to investigate the general abstract initial-boundary value problem
\begin{equation*}\tag{AV}
\left\{\begin{array}{rcll}
\frac{\partial}{\partial t}{u}(t)&=& \Delta u(t),\qquad &t\geq 0,\\
u(t)_{|\Gamma}&\in& {\mathcal Y}, &t\geq 0,\\
\frac{\partial }{\partial t}u(t)_{|\Gamma}&=&P_{\mathcal Y}\left(-\frac{\partial u(t)}{\partial\nu}+\left(\gamma \Delta_\Gamma-{\mathcal S}\right) u(t)_{|\Gamma}\right), &t\geq 0,\\
u(0)&=&u_0,\\
u(0)_{|\Gamma}&=&v_0.
\end{array}
\right.
\end{equation*}
Here $\gamma\in\mathbb R_+$, 
\begin{equation*}\label{mathcals}
{\mathcal S}\in{\mathcal L}\big(H^\frac{1}{2}(\Gamma;H);L^2(\Gamma;H)\big)
\end{equation*}
and $\Delta_\Gamma$ denotes the (dissipative) Laplace--Beltrami operator
on the $(n-1)$-dimensional (differentiable, orientable) manifold $\Gamma$, with the convention that $\gamma=0$ if $n=1$, and hence if $\Gamma$ only consists of isolated points. 
The vector-valued Sobolev space $H^1(\Gamma;H)$ can be defined in the usual way as the vector-valued version of the scalar-valued space $H^1(\Gamma)$ as introduced, i.e., in~\cite[{\S}I.7.3]{LioMag72}. The Laplace operator appearing in $\rm(AV)$ is defined weakly. 

While weak defining the Laplace operator on open domains is standard, a more detailed introduction of the (weakly defined) Laplace--Beltrami operator is in order\footnote{ Observe that any differentiable function $g:\Gamma\to H$ is a mapping between the differentiable manifold $\Gamma$ and the (trivial) Hilbert manifold $H$, whose tangent bundles are $T{\Gamma}\cong {\Gamma}\times {\mathbb R}^{n-1}$ and $TH\cong H\times H$, respectively. Accordingly, at any point $x\in {\Gamma}$ the derivative $\nabla g(x):T_x\Gamma\to T_{g(x)}H$ is a bounded linear operator from ${\mathbb R}^{n-1}$ to $H$ -- hence it can actually be seen as a vector in $H^{n-1}$.}. A definition of the Laplace--Beltrami operator by means of Hilbert space techniques has been performed in the recent preprint~\cite{AreBieEls08}. In fact, 
$$\left(\nabla_\Gamma f(\cdot)|\nabla_\Gamma g(\cdot)\right)_{H^{n-1}}:\Gamma\to{\mathbb C},\qquad f,g\in H^1(\Gamma;H),$$ 
can be defined as the Lebesgue-integrable mapping such that its restriction to any chart $(V,\xi)$ on $\Gamma$ satisfies
$$\left(\nabla_\Gamma f(\cdot)|\nabla_\Gamma g(\cdot)\right)_{{H^{n-1}}_{\big| V}}=\sum_{i,j=1}^{n-1} \left(g^{ij} D_i(f\circ \xi^{-1})\circ \xi\big|D_j(g\circ \xi^{-1})\circ \xi\right)_{H},$$
where $g$ is the canonical Riemannian metric of the surface $\Gamma$. This expression defines in turn a sesquilinear form, and the linear operator associated with this sesquilinear form is the (weakly defined) Laplace--Beltrami operator $\Delta_\Gamma$. We refer to~\cite[\S1]{AreBieEls08} for details.

\begin{rem}
Clearly, both $\Delta$ and $\Delta_\Gamma$ may be replaced by general elliptic operators with real-valued coefficients in pretty much the same way~\cite{MugRom06} generalizes~\cite{AreMetPal03}. Similarly, lower order terms may be added.
\end{rem}

It is known that the right setting for the study of systems of this kind is either the space of continuous functions on $\overline{\Omega}$ or else an $L^p$-product space. We are going to follow the latter approach throughout this note.

\begin{lemma}\label{vdir}
The space 
\begin{equation}\label{VEdindir}
V_{\mathcal Y}:=\left\{{\bf f}:=\begin{pmatrix}f\\ f_{|\Gamma}\end{pmatrix}\in H ^1(\Omega;H)\times \left(H^s(\Gamma;H)\cap{\mathcal Y}\right)\right\}
\end{equation}
is dense in ${\mathcal L}^2_{\mathcal Y}:=L^2(\Omega;H)\times {\mathcal Y}$ for all $s\ge 0$.
\end{lemma}

In no confusion is possible, in the following we will write $\mathcal L^2$ instead of ${\mathcal L}^2_{\mathcal Y}$.

\begin{proof}
This is a slight modification of~\cite[Lemma~5.6]{MugRom06}. More precisely, the assumptions in~\cite[Lemma~5.6]{MugRom06} can be weakened by merely assuming that $H^1(\Gamma;H)\cap {\mathcal Y}$ is \emph{dense} in the range of the trace operator, instead of coinciding with it. This density condition is satisfied by assumption, hence the claim follows.
\end{proof}

In the following we set either $s=1$ if $\gamma>0$, or $s=\frac{1}{2}$ if $\gamma=0$. Accordingly, 
$$V_{\mathcal Y}:=\left\{{\bf f}=\begin{pmatrix}f\\ f_{|\Gamma}\end{pmatrix}\in H ^1(\Omega;H)\times \left(H^1(\Gamma;H)\cap{\mathcal Y}\right)\right\}\qquad\hbox{if }\gamma>0$$
or 
$$V_{\mathcal Y}:=\left\{{\bf f}:=\begin{pmatrix}f\\ f_{|\Gamma}\end{pmatrix}\in H ^1(\Omega;H)\times \left(H^\frac{1}{2}(\Gamma;H)\cap{\mathcal Y}\right)\right\}\qquad\hbox{if } \gamma=0.$$
We consider a form $(a_{\mathcal S},V_{\mathcal Y})$ defined by
\begin{equation*}
 a_{\mathcal S}\left({\bf f},{\bf g}\right):=\int_\Omega \left(\nabla f(x) | \nabla g(x)\right)_{H^n} dx+\gamma\int_{\Gamma} \left (\nabla_\Gamma f(z)|\nabla_\Gamma g(z)\right)_{H^{n-1}} d\sigma(z)+({\mathcal S}f_{|\Gamma}| g_{|\Gamma})_{\mathcal Y},\qquad {\bf f},{\bf g}\in V_{\mathcal Y},
\end{equation*}
where the second addend on the right hand side corresponds to the Laplace--Beltrami operator on the Riemannian manifold $\Gamma$ (recall that by convention $\gamma=0$ whenever $n=1$). We remark that 
$$({\mathcal S}f_{|\Gamma}| g_{|\Gamma})_{\mathcal Y}=({\mathcal S}f_{|\Gamma}| P_{\mathcal Y}g_{|\Gamma})_{\mathcal Y}=( P_{\mathcal Y}{\mathcal S}f_{|\Gamma}| g_{|\Gamma})_{\mathcal Y}\qquad \hbox{for all }{\bf f},{\bf g}\in V_{\mathcal Y},$$
so that the third addend in the definition of $(a_{\mathcal S},V_{\mathcal Y})$ is well-defined.

By a principle presented in~\cite[Appendix]{CarMug09} and based on~\cite[Thm.~4.5.1]{Gra08}, the classical Maz'ya inequality (cf.~\cite[{\S}4.11.2]{Maz85}) can be extended to the vector-valued case. Accordingly, in either case $V_{\mathcal Y}$ is a Hilbert space with respect to the norm defined by
$$({\bf f}|{\bf g})_{V_{\mathcal Y}}:=\int_\Omega(\nabla f(x)|\nabla g(x))_{H^n} dx + \int_{\Gamma} \left(\nabla_\Gamma f(z)|\nabla_\Gamma g(z)\right)_{H^{n-1}}d\sigma(z)+(f_{|\Gamma}|g_{|\Gamma})_{\mathcal Y}\qquad \hbox{if } \gamma>0$$
or
$$({\bf f}|{\bf g})_{V_{\mathcal Y}}:=\int_\Omega(\nabla f(x)|\nabla g(x))_{H^n} dx 
%+ \int_{\Gamma} \left((-\Delta)^\frac{1}{4}u_{|\Gamma}|(-\Delta)^\frac{1}{4}v_{|\Gamma}\right)_{H^{n-1}}d\sigma(z)
+(f_{|\Gamma}|g_{|\Gamma})_{\mathcal Y}\qquad \hbox{if }\gamma=0.$$

\begin{theo}\label{formdir}
The operator $\Delta_{{\mathcal Y},{\mathcal S}}$ associated with $(a_{\mathcal S},V_{\mathcal Y})$ generates an analytic semigroup $(e^{t\Delta_{{\mathcal Y},{\mathcal S}}})_{t\ge 0}$ with angle $\frac{\pi}{2}$ on ${\mathcal L}^2$. 

The operator $\Delta_{{\mathcal Y},{\mathcal S}}$ is dissipative if the operator ${\mathcal S}$ is accretive and in this case the semigroup $(e^{t\Delta_{{\mathcal Y},{\mathcal S}}})_{t\ge 0}$ is contractive. The operator $\Delta_{{\mathcal Y},{\mathcal S}}$  is self-adjoint if and only if the operator ${\mathcal S}$ is self-adjoint and in this case the semigroup $(e^{t\Delta_{{\mathcal Y},{\mathcal S}}})_{t\ge 0}$ is self-adjoint. The operator $\Delta_{{\mathcal Y},{\mathcal S}}$ has compact resolvent if and only if $\Omega,\Gamma$ have finite measure, provided that $H$ is finite dimensional: in this case the semigroup $(e^{t\Delta_{{\mathcal Y},{\mathcal S}}})_{t\ge 0}$ is compact.
\end{theo}

The proof is based on the approach presented, e.g., in~\cite[Chapt.~VI]{DauLio88}. We borrow our terminology from~\cite{Are06}.

\begin{proof}
We are going to show that $(a_{\mathcal S},V_{\mathcal Y})$ is associated with an operator that generates a cosine family with phase space $V_{\mathcal Y}\times {\mathcal L}^2$ in the sense of~\cite[{\S}3.14]{AreBatHie01}. To this aim, we show that for all $\gamma\in {\mathbb R}_+$ the densely defined sesquilinear form $(a_{\mathcal S},V_{\mathcal Y})$ is continuous and elliptic (with respect to ${\mathcal L}^2$), i.e.,
$${\rm Re}a_{\mathcal S}({\bf f},{\bf f}) +\omega \|{\bf f}\|^2_{{\mathcal L}^2} \ge \alpha \|{\bf f}\|^2_{V_{\mathcal Y}} \qquad\hbox{for all }{\bf f}\in V_{\mathcal Y}$$
for some $\alpha>0$ and a suitable $\omega\in\mathbb R$.

Continuity follows from the Cauchy--Schwarz inequality. Ellipticity (with respect to ${\mathcal L}^2$) follows from ellipticity (with respect to $L^2(\Omega;H)$ and $L^2(\Gamma;H)$) of the forms associated with the Laplace and Laplace--Beltrami operators, corresponding to the first two addends of $(a_{\mathcal S},V_{\mathcal Y})$. The third addend in the definition of $a_{\mathcal S}$ is sesquilinear and defined on $H^\frac{1}{2}(\Gamma;H)\times H^\frac{1}{2}(\Gamma;H)$, hence it can be neglected by a perturbation argument (see~\cite[Lemma~2.1]{Mug08}). Finally, because 
\begin{eqnarray*}
|{\rm Im} a_{\mathcal S}({\bf f},{\bf f})|&=&|{\rm Im}({\mathcal S}f_{|\Gamma}|f_{|\Gamma})_{\mathcal Y}|\le \|S\| \|f\|_{H^\frac{1}{2}(\Gamma;H)}\|f_{|\Gamma}\|_{L^2(\Gamma;H)}\\ 
&\le& M\|S\| \|f\|_{H^1(\Omega;H)}\|f_{|\Gamma}\|_{L^2(\Gamma;H)}, 
\end{eqnarray*}
for some $M>0$ and all ${\bf f}\in V_{\mathcal Y}$ due to boundedness of the trace operator from $H^1(\Omega;H)$ to $H^\frac{1}{2}(\Gamma;H)$, the announced generation of a cosine family follows by~\cite[Thm.~4]{Cro07}. It is known that generators of cosine operator functions also generate analytic semigroups with angle $\frac{\pi}{2}$, see~\cite[Thm.~3.14.17]{AreBatHie01}.

Because the forms associated with the Laplace and Laplace--Beltrami operators are accretive, accretivity of $(a_{\mathcal S},V_{\mathcal Y})$ is clear provided ${\mathcal S}$ is accretive. A direct computation shows that $(a_{\mathcal S},V_{\mathcal Y})$ is symmetric if and only if ${\mathcal S}$ is self-adjoint. The assertion on compactness follows from the Aubin--Lions Lemma, see~\cite[Prop.~III.1.3]{Sho97}.
\end{proof}

The proof of the following is based on~\cite[Rem.~2.2]{AreMetPal03}. 

\begin{prop}\label{identdir}
Assume $\Omega$ to have $C^2$-boundary. For all $\gamma\in {\mathbb R}_+$ and ${\mathcal S}\in{\mathcal L}(L^2(\Gamma;H))$ the operator $\Delta_{{\mathcal Y},{\mathcal S}}$ associated with
$(a_{\mathcal S},V_{\mathcal Y})$ is given by
\begin{eqnarray*}
D(\Delta_{{\mathcal Y},{\mathcal S}})&=&\left\{{\bf f}:=\begin{pmatrix} f\\ f_{|\Gamma}\end{pmatrix}\in V_{\mathcal Y}: \Delta f\in L^2(\Omega;H),\; \Delta_\Gamma f_{|\Gamma}\in L^2(\Gamma;H),\hbox{ and }\frac{\partial f}{\partial \nu}\in L^2(\Gamma;H)\right\},\\
\Delta_{{\mathcal Y},{\mathcal S}}&=&\begin{pmatrix} \Delta & 0\\
-P_{\mathcal Y}\frac{\partial}{\partial \nu} & P_{\mathcal Y}\left(\gamma\Delta_\Gamma -{\mathcal S}\right)\end{pmatrix},
\end{eqnarray*}
hence $(e^{t\Delta_{{\mathcal Y},{\mathcal S}}})_{t\ge 0}$ yields the solution to ${\rm(AV)}$. If in particular ${\bf f}\in D(\Delta_{{\mathcal Y},{\mathcal S}})$, then $f\in H^\frac{3}{2}(\Omega;H)\cap H^2_{loc}(\Omega;H)$.
\end{prop}

Observe that in general ${A}_{\mathcal Y}$ would not operate on ${\mathcal L}^2$ if we would drop the term $P_{\mathcal Y}$.

\begin{proof}
By definition, the operator associated with $(a_{\mathcal S},V_{\mathcal Y})$ is given by
\begin{eqnarray*}
D(B_{{\mathcal Y},{\mathcal S}})&:=&\left\{{\bf f}\in V_{\mathcal Y} :\exists {\bf g}
\in {\mathcal L}^2 \hbox{ s.t. } {a}\left({\bf f},{\bf h}\right)=\left({\bf g}| {\bf h}\right)_{{\mathcal L}^2}\; \forall {\bf h}\in V_{\mathcal Y}\right\},\\
B_{{\mathcal Y},{\mathcal S}}{\bf f}&:=&-{\bf g}.
\end{eqnarray*}

By the perturbation theorem of Desch--Schappacher (cf.~\cite{DesSch84}) a relatively bounded perturbation does not affect the domain of an operator. Hence we can assume w.l.o.g.\ that 
$${\mathcal S}=0.$$ 
In order to prove that $\Delta_{{\mathcal Y},{\mathcal S}}\subset B_{{\mathcal Y},{\mathcal S}}$ take
${\bf f},{\bf h}\in V_{\mathcal Y}$. By the Gau{\ss}--Green formulae and the (weak) definition of the Laplace--Beltrami operator we obtain
\begin{align*}
a_{\mathcal S}\left({\bf f},{\bf h}\right)
&=\int_{\Omega} (\nabla f(x)|\nabla h(x))_{H^n} dx+\gamma\int_{\Gamma} \left (\nabla_\Gamma f_{|\Gamma}(z)|\nabla_\Gamma h_{|\Gamma}(z)\right)_{H^{n-1}} d\sigma(z)\\
&= -\int_{\Omega} (\Delta f(x)| h(x))_H dx\\
&\qquad + \int_{\Gamma} \left(\frac{\partial f}{\partial \nu}f(z)|h_{|\Gamma}(z)\right)_H d\sigma(z) -\gamma\int_{\Gamma} \left (\Delta_{\Gamma} f_{|\Gamma}(z)| g_{|\Gamma}(z)\right)_H d\sigma(z)\\
&= -\int_{\Omega} (\Delta f(x)| h(x))_H dx+ \left(\frac{\partial f}{\partial \nu}f(z)|h(z)\right)_{{\mathcal Y}} -\gamma \left (\Delta_\Gamma f_{|\Gamma}| g_{|\Gamma}\right)_{{\mathcal Y}}\\
&= \left( \begin{pmatrix}-\Delta f\\
P_{\mathcal Y}\left(\frac{\partial f}{\partial\nu}-\gamma \Delta_\Gamma f_{|\Gamma}\right) \end{pmatrix}\big|
\begin{pmatrix}h\\ h|_{{\Gamma}}\end{pmatrix}\right)_{{\mathcal L}^2}=:\left({\bf g}| {\bf h}\right)_{{\mathcal L}^2},
\end{align*}
and the operator $\Delta_{{\mathcal Y},{\mathcal S}}$ has the claimed form.

Conversely, let ${\bf f}\in D(B_{{\mathcal Y},{\mathcal S}})$. The above computation also shows that $\Delta f$ and $\Delta_\Gamma f_{|\Gamma}$ are well defined elements of $L^2(\Omega;H)$ and $L^2(\Gamma;H)$, respectively, and that $f$ has a weak normal derivative in $L^2(\Gamma;H)$. We deduce that $f\in H ^\frac{3}{2}(\Omega;H)$ by~\cite[Thm.~2.7.4]{LioMag72} -- suitably extended to the vector-valued case by virtue of~\cite[Thm.~4.5.1]{Gra08}.
\end{proof}

\begin{rem}
The vectors in $D(A^2_{\mathcal Y})$ also satisfy the additional boundary condition
\begin{equation}\label{lapbel}
(\Delta u)_{|\Gamma}\in {\mathcal Y}\quad\hbox{and}\quad
(\Delta u)_{|\Gamma}+P_{\mathcal Y}\frac{\partial u}{\partial \nu}+P_{\mathcal Y}\left({\mathcal S}u_{|\Gamma}-\gamma \Delta_\Gamma u_{|\Gamma}\right)=0 
\end{equation}
for all $z\in\Gamma$. Conditions~\ref{lapbel} can be interpreted as a formulation of Wentzell--Robin boundary conditions which is stronger than the dynamic one that is usual in the context of $L^p$-spaces. Due to the regularising effect of the analytic semigroup $(e^{t\Delta_{{\mathcal Y},{\mathcal S}}})_{t\ge 0}$, these additional conditions are satisfied by the solution $\rm(AV)$ for any time $t>0$.
% 2) The variational structure of the problem considered in this paper has great advantages. In particular, via the Gau{\ss}--Green formula we can make use of Hilbert space techniques. If however a multiplicative coefficient $\rho\in\mathbb R$ is considered in the lower-left entry of $\Delta_{{\mathcal Y},{\mathcal S}}$, then the variational structure is lost. Does this modified operator generates a semigroup, too? In the scalar case it is known that the answer is yes if and only if $n=1$ or $\rho>0$, cf.~\cite{Mug08c,VazVit08} and~\cite[{\S}15]{Bel94}. It can be conjectured that this result holds true in the vector-valued case. 
\end{rem}

Consider a sequence $({\mathcal Y}_n)_{n\in\mathbb N}$ of closed subspaces of $L^2(\Gamma;H)$ such that the associated sequence of orthogonal projections $(P_{{\mathcal Y}_n})_{n\in\mathbb N}$ converges in operator norm. Then its limit is also necessarily a projection and a contraction, i.e., an orthogonal projection -- say, onto a subspace $\mathcal Y$. 
Consider moreover a sequence $({\mathcal S}_n)_{n\in\mathbb N}$ in ${\mathcal L}(H^\frac{1}{2}(\Gamma;H);L^2(\Gamma;H))$ that converges in operator norm to some ${\mathcal S}\in {\mathcal L}(H^\frac{1}{2}(\Gamma;H);L^2(\Gamma;H))$.  
Now, it is quite natural to conjecture that \textit{$\Delta_{{\mathcal Y}_n,{\mathcal S}_n}$ converges to $\Delta_{{\mathcal Y},{\mathcal S}}$ in a suitable sense}. 

Observe that no kind of convergence from above or below of the form family $(a_{{\mathcal S}_n},V_{{\mathcal Y}_n})_{n\in\mathbb N}$ holds -- in our case one typically has $V_{{\mathcal Y}_n}\cap V_{{\mathcal Y}_m}=V_{{\mathcal Y}_0}$ for some lower-dimensional ${\mathcal Y}_0$, whenever $n\not=m$ -- so that in general $V_{{\mathcal Y}_0}$ is not dense in any $V_{{\mathcal Y}_n}$. Furthermore, the operators $\Delta_{{\mathcal Y}_n,{\mathcal S}_n}$ and $\Delta_{{\mathcal Y},{\mathcal S}}$ act on $L^2(\Omega;H)\times {\mathcal Y}_n$ and $L^2(\Omega;H)\times {\mathcal Y}$, respectively, i.e., they generally act on different spaces. All in all, it seems that well-known results for convergence of operators associated with forms (e.g., those due to Kato and Simon) cannot be applied to our setting. Some results on approximation of operators acting on different spaces have been recently obtained by Ito and Kappel (see e.g.~\cite[Chapt.~4]{ItoKap02}), but it seems that they fall short of fitting our framework, too.

The different approach proposed by Post in~\cite{Pos06} and further developed in~\cite{MugNitPos10} seems to be more appropriate. In order to apply Post's results, we need to impose a structural assumption on $\mathcal Y$ that will prove a significant simplification in our framework.

\begin{theo}\label{olafrev}
Let $(Y_n)_{n\in\mathbb N}$ be a sequence of closed subspaces of $H$. Consider a further closed subspace $Y$ of $H$ and a family $({J}^{\downarrow n})_{n\in\mathbb N}$ of unitary operators on $H$ such that ${J}^{\downarrow n}Y_n=Y$ for all $n\in\mathbb N$. Assume furthermore that $\lim_{n\to\infty} {J}^{\downarrow n}={\rm Id}$ in operator norm and consider the spaces 
\begin{eqnarray}
{\mathcal Y}&:=&\{f\in L^2(\Gamma;H): f(z)\in Y \hbox{ for a.e. } z\in\Gamma\},\qquad \hbox{and}\label{oneY0}\\
{\mathcal Y}_n&:=&\{f\in L^2(\Gamma;H): f(z)\in Y_n \hbox{ for a.e. } z\in\Gamma\},\qquad n\in\mathbb{N}.\label{oneY0b}
\end{eqnarray}
Let $(S_n)_{n\in\mathbb N}$ be a sequence of accretive bounded linear operators on $H$ that converges in operator norm to some $S\in{\mathcal L}(H)$ and define linear operators ${\mathcal S}_n,S\in {\mathcal L}(H^\frac{1}{2}(\Gamma;H))$ by
$${\mathcal S}_n g:= S_n\circ g,\; n\in\mathbb N,\qquad\hbox{and}\qquad {\mathcal S}g:= S\circ g,\qquad g\in H^\frac{1}{2}(\Gamma;H).$$
Then both families $(R(\lambda,\Delta_{{\mathcal Y}_n,{\mathcal S}_n}))_{n\in\mathbb N}$ and $(e^{t\Delta_{{\mathcal Y}_n,{\mathcal S}_n}})_{n\in\mathbb N}$ of bounded linear operators on ${\mathcal L}^2_{{\mathcal Y}_n}$ converge in operator norm to the bounded linear operators $R(\lambda,\Delta_{{\mathcal Y},{\mathcal S}})$ and to $e^{t\Delta_{{\mathcal Y},{\mathcal S}}}$ on ${\mathcal L}^2_{\mathcal Y}$, for all ${\rm Re}\lambda>0$ and for all $t>0$ respectively. Moreover, if $H$ is finite dimensional and $\Omega,\Gamma$ have finite measure, then the (discrete) spectrum of $\Delta_{{\mathcal Y}_n,{\mathcal S}_n}$ converges to the (discrete) spectrum of $\Delta_{{\mathcal Y},{\mathcal S}}$.
\end{theo}

\begin{rem}
Observe that the phenomenon observed in Theorem~\ref{olafrev} is intrinsically related to the vector-valued case. If in fact ${\rm dim}\; H=1$, then each sequence $(Y_n)_{n\in\mathbb N}$ of subspaces of $H={\mathbb C}$ such that $(P_{Y_n})_{n\in\mathbb N}$ converges is eventually constant -- with value either $\{0\}$ or $H$ -- so that the assertion becomes trivial.
\end{rem}

The proof is based on an abstract convergence scheme discussed in~\cite[Appendix]{Pos06}, which we briefly recall for the sake of self-containedness. The following collects results from~\cite[Thms.~A.5 and~A.10]{Pos06}

\begin{prop}\label{postlemma}
Let $\mathcal H,\mathcal H_1,\tilde{\mathcal H},\tilde{\mathcal H_1}$ be Hilbert spaces such that $\mathcal H_1\hookrightarrow \mathcal H$ and $\tilde{\mathcal H}_1\hookrightarrow \tilde{\mathcal H}$ with dense embeddings. Let ${\mathfrak h}:{\mathcal H}_1\times {\mathcal H}_1\to{\mathbb C}$ and $\tilde{\mathfrak h}:\tilde{\mathcal H}_1\times \tilde{\mathcal H}_1\to{\mathbb C}$ be continuous, accretive and elliptic (with respect to $\mathcal H$ and $\tilde{\mathcal H}$, respectively) with associated operators $\mathfrak A$ and $\tilde{\mathfrak A}$. Consider operators ${\mathcal J}\in{\mathcal L}({\mathcal H},\tilde{\mathcal H})$, $\tilde{\mathcal J}\in{\mathcal L}(\tilde{\mathcal H},{\mathcal H})$, ${\mathcal J}_1\in{\mathcal L}({\mathcal H}_1,\tilde{\mathcal H}_1)$, $\tilde{\mathcal J}_1\in{\mathcal L}(\tilde{\mathcal H}_1,{\mathcal H}_1)$. Let moreover the above spaces and operators satisfy the following conditions:
\begin{eqnarray}
\| {\mathcal J} {\bf f}-{\mathcal J}_1 {\bf f}\|_{\tilde{\mathcal H}} &\le& \delta \|{\bf f}\|_{{\mathcal H}_1},\label{post1}\\ 
\| \tilde{\mathcal J} {\bf u}-\tilde{\mathcal J}_1 {\bf u}\|_{{\mathcal H}} &\le& \delta \|{\bf u}\|_{\tilde{\mathcal H}_1},\label{post2}\\
| ({\mathcal J}{\bf f}| {\bf u})_{\tilde{\mathcal H}}-({\bf f}| \tilde{\mathcal J}{\bf u})_{{\mathcal H}}| &\le& \delta \|{\bf f}\|_{{\mathcal H}} 
\|{\bf u}\|_{\tilde{\mathcal H}},\label{post3}\\
| \tilde{\mathfrak h}({\mathcal J}_1{\bf f}| {\bf u})-{\mathfrak h}({\bf f}| \tilde{\mathcal J}_1{\bf u})| &\le& \delta \|{\bf f}\|_{{\mathcal H}_1} \|{\bf u}\|_{\tilde{\mathcal H}_1},\label{post4}\\
\| {\bf f}-\tilde{\mathcal J}{\mathcal J}{\bf f}\|_{\mathcal H}&\le& \delta \|{\bf f}\|_{{\mathcal H}_1},\label{post5}\\
\| {\bf u}-{\mathcal J}\tilde{\mathcal J}{\bf u}\|_{\tilde{\mathcal H}}&\le& \delta \|{\bf u}\|_{\tilde{\mathcal H}_1},\label{post6}\\
\|{\mathcal J}{\bf f}\|_{\tilde{\mathcal H}}&\le& 2\|{\bf f}\|_{\mathcal H},\label{post7}\\
\|\tilde{\mathcal J}{\bf u}\|_{{\mathcal H}}&\le& 2\|{\bf u}\|_{\tilde{\mathcal H}}.\label{post8}
\end{eqnarray}
for some $\delta>0$. Then
$$\|R(\lambda,\tilde{{\mathfrak A}})-{\mathcal J}R(\lambda,{{\mathfrak A}})\tilde{\mathcal J}\|\le M\delta$$
for some $M>0$.
\end{prop}

We emphasize that the convergence assertion is rather poor at a numerical level, but fairly strong at a functional analytical level: it states convergence in operator norm, rather than just strong convergence as done e.g.\ by the various Trotter--Kato-type theorems. We are now in the position to prove Theorem~\ref{olafrev}.

\bigskip

\noindent
\emph{Proof of Theorem~\ref{olafrev}.} Fix $n\in\mathbb N$. We apply Proposition~\ref{postlemma} setting 
$$\mathcal H:={\mathcal L}^2_{{\mathcal Y}_n},\;\tilde{\mathcal H}:={\mathcal L}^2_{\mathcal Y},\;\mathcal H_1:=V_{{\mathcal Y}_n},\;\tilde{\mathcal H}_1:=V_{\mathcal Y},$$
along with
$${\mathfrak h}:=(a_{{\mathcal S}_n}, V_{{\mathcal Y}_n})\hbox{ and }\tilde{\mathfrak h}:=(a_{{\mathcal S}},V_{{\mathcal Y}}).$$
Observe that accretivity of ${\mathfrak h},\tilde{\mathfrak h}$ follows from accretivity of the operators $S_n,S$.
Define moreover ${\mathcal J}\in{\mathcal L}({\mathcal H},\tilde{\mathcal H})$ by
\begin{equation}\label{j1}
{\mathcal J}{\bf f}:=\begin{pmatrix}
J^{\downarrow n}\circ f_1\\
J^{\downarrow n}\circ f_2\\
    \end{pmatrix},\qquad {\bf f}:=\begin{pmatrix}f_1\\ f_2\end{pmatrix}\in {\mathcal H},
\end{equation}
and
\begin{equation}\label{j2}
\tilde{\mathcal J}{\bf u}:={\mathcal J}^{-1}{\bf u}=\begin{pmatrix}
(J^{\downarrow n})^{-1}\circ u_1\\
(J^{\downarrow n})^{-1}\circ u_2\\
    \end{pmatrix},\qquad {\bf u}:=\begin{pmatrix}u_1\\ u_2\end{pmatrix}\in \tilde{\mathcal H},
\end{equation}
and moreover ${\mathcal J}_1:={\mathcal J}$ and $\tilde{\mathcal J}_1:=\tilde{\mathcal J}$. 

It is apparent that~\eqref{post1}, \eqref{post2}, \eqref{post7} and \eqref{post8} are trivially satisfied for $\delta$ large enough (and getting smaller and smaller as $n$ increases). Moreover,~\eqref{post3},\eqref{post5} and \eqref{post6} hold because ${\mathcal J},\tilde{\mathcal J}$ are unitary with ${\mathcal J}^*=\tilde{\mathcal J}$. 

Finally, observe that by~\eqref{oneY0},~\eqref{oneY0b},~\eqref{j1},~\eqref{j2} the operators ${\mathcal J}_1$, $\tilde{\mathcal J}_1$ do not depend on space, hence they commute with the local operators associated with the forms ${\mathfrak h},\tilde{\mathfrak h}$. Furthermore, for all ${\bf f}\in{\mathcal H}_1$ and all ${\bf u}\in\tilde{\mathcal H}_1$
$$({\mathcal S} {\mathcal J}^{\downarrow n} f_{|\Gamma}|u_{|\Gamma})_{\mathcal Y} - ({\mathcal S}_n f_{|\Gamma}|({\mathcal J}^{\downarrow n})^{-1} u_{|\Gamma})_{{\mathcal Y}_n}=(({\mathcal S} {\mathcal J}^{\downarrow n}-{\mathcal J}^{\downarrow n}{\mathcal S}_n) f_{|\Gamma}| u_{|\Gamma})_{{\mathcal Y}_n}$$
which converges to $0$ because ${\mathcal J}^{\downarrow n\; -1}{\mathcal S}_n{\mathcal J}^{\downarrow n}$ converges to $\mathcal S$ in operator norm. We conclude that~\eqref{post4} is satisfied.

Then the convergence of $(R(\lambda,\Delta_{\mathcal Y_n}))_{n\in\mathbb N}$ follows from Proposition~\ref{postlemma}. The remaining assertions follow from~\cite[Thms.~A.10 and A.11]{Pos06}.
\hfill\qed

\begin{rem}
Let us consider the case of a more general diffusion equation of the form
\begin{equation}\tag{AV$_D$}
\left\{\begin{array}{rcll}
\frac{\partial}{\partial t}{u}(t)&=& \nabla\cdot(D\nabla u(t)),\qquad &t\geq 0,\\
u(t)_{|\Gamma}&\in& {\mathcal Y}, &t\geq 0,\\
\frac{\partial }{\partial t}u(t)_{|\Gamma}&=&
-P_{\mathcal Y}\frac{\partial_D u(t)}{\partial\nu}+\left(\gamma\Delta_\Gamma-{\mathcal S})\right) u(t)_{|\Gamma}, &t\geq 0,\\
u(0)&=&u_0,\\
u(0)_{|\Gamma}&=&v_0,
\end{array}
\right.
\end{equation}
where $D\in C^1(\overline{\Omega};{\mathcal L}(H^n))$ satisfies for some $\mu>0$ the following ellipticity condition: 
$${\rm Re}(D(x)\xi| \xi)_{H^n}\geq \mu \| \xi\|_{H^n}^2\qquad\hbox{ for all } \xi\in H^n\hbox{ and all\ }x\in\overline{\Omega}.$$	
The subspace $\mathcal Y$ as well as the operator $\mathcal S$ are now fixed. Then, a variational approach can still be pursued, after introducing suitable weighted Bochner spaces ${\mathcal L}^2_{\mathcal Y,D}$ as it has been done in~\cite{MugRom06}. Due to uniform ellipticity, the coefficients do not degenerate on the boundary, yielding that $\|\cdot\|_{{\mathcal L}^2_{\mathcal Y,D}}$ and $\|\cdot\|_{{\mathcal L}^2_{\mathcal Y}}$ are equivalent norms on $\mathcal L^2_{\mathcal Y,D}$.

Now, consider a uniformly elliptic family $(D_k)_{k\in\mathbb N}\subset C^1(\overline{\Omega};{\mathcal L}(H^n))$ of coefficients such that $D_k(x)$ is self-adjoint for all $x\in\Omega$ and all $k\in\mathbb N$. Consider the sesquilinear form $a_k$ arising from the problem $({\rm AV}_{D_k})$, $k\in\mathbb N$, whose domains all coincide with $V_{\mathcal Y}$. Denote by $\Delta_k$ the associated operator. These operators are uniformly sectorial -- actually, all their numerical ranges are contained in the negative halfline. If the sequence $(D_k)_{k\in\mathbb N}$ converges strongly, then $(a_k({\bf f},{\bf f}))_{k\in\mathbb N}$ is a Cauchy sequence for all ${\bf f}\in V_{\mathcal Y}$. Therefore, by a known result due to Kato (see~\cite[{\S}VIII.3]{Kat66}), $(R(\lambda,\Delta_{\mathcal Y_k}))_{k\in\mathbb N}$ converges strongly for all ${\rm Re}\lambda>0$. By simple functional calculus arguments this also implies strong convergence of $(e^{z\Delta_k})_{k\in\mathbb N}$ for all $z$ in the open right halfplane. This is comparable with~\cite[Thm.~3.1]{CocFavGol08}. A similar assertion concerning convergence in operator norm can also be obtained applying Proposition~\ref{postlemma}.
\end{rem}

\section{Lattice-based invariance properties}\label{posdom}

This section is devoted to the characterisation of qualitative properties of $(e^{t\Delta_{{\mathcal Y},{\mathcal S}}})_{t\ge 0}$. These can often be discussed in terms of invariance of relevant subsets of the state space $\mathcal L^2$ -- most notably, order intervals\footnote{ It has been observed in~\cite[{\S}5]{CarMugNit08} that also invariances of some subspaces of the state space often reveal important properties of the evolution equation. In fact, all results in this section also apply when order intervals are replaced by subspace -- of course, even dropping all lattice assumptions.}. By Ouhabaz's well-known invariance criterion, such invariance properties can be characterized by simple, almost linear algebraic properties of a quadratic form. In a more general form presented in~\cite[Thm.~2.1]{ManVogVoi05}, Ouhabaz's criterion can be stated as follows.

\begin{lemma}\label{ouhalemma}
Let $\mathcal H$ be a separable Hilbert space and $\bf a$ a sesquilinear form  with dense domain $\mathcal V$ that is continuous and elliptic with respect to $\mathcal H$. A closed convex set $\mathcal C$ of $\mathcal H$ is invariant under the semigroup associated with $a$ if and only if $\mathcal V$ is invariant under the orthogonal projection $\mathcal P$ of $\mathcal H$ onto $\mathcal C$ and moreover ${\rm Re}\;{\bf a}({\mathcal P}{\mathfrak u},{\mathfrak u}-{\mathcal P}{\mathfrak u})\ge 0$ for all ${\mathfrak u}\in\mathcal V$.
\end{lemma}

In the remainder of this section assume for simplicity that $\gamma>0$, i.e.,
$$
V_{\mathcal Y}=\left\{{\bf f}\in H ^1(\Omega;H)\times H^1(\Gamma;H):f_{|\Gamma}\in{\mathcal Y}\right\}.
$$
(Still, all assertions hold true in the case $\gamma=0$ with obvious, minor modifications in the proofs).

\bigskip
To warm up, we start by characterising reality of $(e^{t\Delta_{{\mathcal Y},{\mathcal S}}})_{t\geq 0}$. A function in ${\mathcal L}^2$ is called $H_{\mathbb R}$-valued if it takes values in the real Hilbert space $H_{\mathbb R}$ underlying $H$ for a.e.\ $x\in\Omega\oplus \Gamma$. 
As a direct consequence of locality the forms associated with the Laplace and Laplace--Beltrami operators we obtain the following.

\begin{prop}
The semigroup $(e^{t\Delta_{{\mathcal Y},{\mathcal S}}})_{t\geq 0}$ leaves invariant the real part of ${\mathcal L}^2$, i.e., the set of all $H_{\mathbb R}$-valued functions in ${\mathcal L}^2$, if and only if $({\mathcal S}{\rm Re}g|{\rm Im}g)_{\mathcal Y}\in\mathbb R$ for all $g\in \mathcal Y$, hence if and only if 
$${\mathcal S} \{f\in H^\frac{1}{2}(\Gamma;H):f(x)\in {H}_{\mathbb R} \hbox{ for a.e. } x\in\Gamma\}\subset \{f\in L^2(\Gamma;H):f(x)\in {H}_{\mathbb R} \hbox{ for a.e. } x\in\Gamma\}.$$
\end{prop}

In typical applications the space $H$ is a Hilbert lattice\footnote{ For the necessary notions from the theory of Banach lattices we refer to~\cite{Nag86,Mey91}. Consequently, also $L^2(\Omega;H)$, $\mathcal Y$ and ${\mathcal L}^2$ are Hilbert lattices. Whenever we refer to an operator on a Hilbert lattice as ``positive'', we always mean ``positivity preserving''.} -- hence we will assume henceforth that 
$$H\cong L^2(\Xi;{\mathbb C})$$ 
for a suitable finite measure space $\Xi$, cf.~\cite[Cor.~2.7.5]{Mey91}. In particular, the scalar products of $L^2(\Omega;H)$ and $L^2(\Gamma;H)$ read now
$$(f|g)_{L^2(\Omega;H)}:=\int_\Omega \int_\Xi f(x,\xi)\overline{g(x,\xi)} d\xi dx\quad\hbox{and}\quad(f|g)_{L^2(\Gamma;H)}:=\int_\Gamma \int_\Xi f(x,\xi)\overline{g(x,\xi)} d\xi dx,$$
respectively. We can define the positive and negative parts and the absolute value of functions in $L^2(\Omega;H)$ pointwise, exploiting the lattice structure of $H$: if 
$u\in L^2(\Omega;H)\cong L^2(\Omega\times \Xi;{\mathbb C})$, then $u^+,u^-$ are the function
$$\Omega\ni x\mapsto u^+(x,\cdot)\in L^2(\Xi;{\mathbb C})\quad\hbox{and}\quad \Omega\ni\omega\mapsto u^-(\omega)\in L^2(\Xi;{\mathbb C}),$$
respectively, where $u^+(x,\cdot),u^-(x,\cdot)$ are well-defined elements of $H=L^2(\Xi,{\mathbb C})$ because $u(x,\cdot)\in L^2(\Xi;{\mathbb C})$. Observe that the orthogonal projection ${\mathcal P}_+$ onto the positive cone of $L^2(\Omega;H)$ acts on any $u\in L^2(\Omega;H)$ as the composition $P_+\circ u$, where $P_+$ is the orthogonal projection onto the positive cone of $H$. 
%This turns out to be relevant when it comes to applying abstract results on differentiation of Lipschitz mappings between Hilbert spaces (like the projection onto the positive cone of $H$), cf.~\cite{Aro76} and references therein.

Let $a,b\in L^2(\Omega; H)\cong L^2(\Omega\times \Xi;{\mathbb C})$ and consider the unbounded order intervals
\begin{eqnarray*}
%[a,b]_{L^2(\Omega;H)}&:=&\left\{f \in L^2(\Omega;H): \min\{a(x),b(x)\}\le f(x)\le \max\{a(x),b(x)\}\hbox{ for a.e. } x\in \Omega\right\},\\
\; [a,+\infty)_{L^2(\Omega;H)}&:=&\left\{f \in L^2(\Omega;H): a(x)\le f(x)\hbox{ for a.e. } x\in \Omega\right\}\\
&\cong &\left\{f \in L^2(\Omega\times\Xi;{\mathbb C}): a(x,\xi)\le f(x,\xi)\hbox{ for a.e. } (x,\xi)\in \Omega\times\Xi\right\},\\
(-\infty,b]_{L^2(\Omega;H)}&:=&\left\{f \in L^2(\Omega;H): f(x)\le b(x)\hbox{ for a.e. } x\in \Omega\right\}\\
&\cong &\left\{f \in L^2(\Omega\times\Xi;{\mathbb C}): f(x,\xi)\le b(x,\xi)\hbox{ for a.e. } (x,\xi)\in \Omega\times\Xi\right\}.
\end{eqnarray*}
These subsets of $L^2(\Omega;H)$ are closed and convex. Similarly, for $c,d\in L^2(\Gamma; H)\cong L^2(\Omega\times \Xi;{\mathbb C})$ one considers the unbounded order intervals $ [c,+\infty)_{\mathcal Y}, (-\infty,d]_{\mathcal Y}$.

\begin{lemma}\label{gteinelemma}
Let $u,v\in H^1(\Omega;H)$ and $\tilde{u},\tilde{v}\in H^1(\Gamma;H)$. Then also $\max\{u,v\}\in H^1(\Omega;H)$ as well as $\max\{\tilde{u},\tilde{v}\}\in H^1(\Gamma;H)$. Furthermore,
\begin{eqnarray}
\nabla \max\{u,v\}(\cdot,\xi)&=&{\mathbf 1}_{\{u(\cdot,\xi)\ge v(\cdot,\xi)\}}\nabla u(\cdot,\xi)+\mathbf 1_{\{u(\cdot,\xi)< v(\cdot,\xi)\}}\nabla v(\cdot,\xi),\label{gteine}\\
\nabla (u-v)^-(\cdot,\xi)&=&\mathbf 1_{\{u(\cdot,\xi)< v(\cdot,\xi)\}}\nabla v(\cdot,\xi),\label{gteine2}\\
\nabla \max\{\tilde{u},\tilde{v}\}(\cdot,\xi)&=&{\mathbf 1}_{\{\tilde{u}(\cdot,\xi)\ge \tilde{v}(\cdot,\xi)\}}\nabla \tilde{u}(\cdot,\xi)+\mathbf 1_{\{\tilde{u}(\cdot,\xi)< \tilde{v}(\cdot,\xi)\}}\nabla \tilde{v}(\cdot,\xi),\label{gteine3}\\
\nabla (\tilde{u}-\tilde{v})^-(\cdot,\xi)&=&\mathbf 1_{\{\tilde{u}(\cdot,\xi)< \tilde{v}(\cdot,\xi)\}}\nabla \tilde{v}(\cdot,\xi),\label{gteine4}
\end{eqnarray}
for a.e.\ $\xi\in\Xi$.
\end{lemma}

Observe that~\eqref{gteine},\eqref{gteine2},\eqref{gteine3},\eqref{gteine4} represent equalities of functions in $L^2(\Omega;{\mathbb C})$ and $L^2(\Gamma;{\mathbb C})$, respectively. In particular, each ``slice'' $u(\cdot,\xi)$ defines a scalar-valued function on $\Omega$: it is the differential of this slice-function that is denoted by $\nabla u(\cdot,\xi)$. The same is valid for $v,\tilde{u},\tilde{v}$.

\begin{proof}
The proof goes in several steps. We will repeatedly use the fact that 
$$u\in L^2(\Omega;H)\cong L^2(\Omega;{\mathbb C})\otimes H\cong L^2(\Omega;{\mathbb C})\otimes L^2(\Xi;{\mathbb C})\cong L^2(\Omega\times\Xi;{\mathbb C})$$
in order to reduce a vector-valued relation to a collection of scalar-valued ones: this follows from the elementary theory of Hilbert tensor products.

1) First of all, we recall the following vector-valued extension of~\cite[Prop.~IX.3]{Bre83}, observed in~\cite[Appendix A]{CarMug09}: \emph{Let $G: H \to H$ be a  Lipschitz continuous mapping and $f \in H^1(\Omega;H)$. If $G(0)=0$, then $G\circ f \in H^1(\Omega;H)$.} The proof is an easy modification of~\cite[Prop.~IX.3]{Bre83}.
In particular, this result applies to the case where $G$ is an orthogonal projection onto an order interval
\begin{eqnarray*}
-%[a,b]_{L^2(\Omega;H)}&:=&\left\{f \in L^2(\Omega;H): \min\{a(x),b(x)\}\le f(x)\le \max\{a(x),b(x)\}\hbox{ for a.e. } x\in \Omega\right\},\\
\; [\alpha,+\infty)_{L^2(\Omega;H)}&:=&\left\{f \in L^2(\Omega;H): \alpha\le f(x)\hbox{ for a.e. } x\in \Omega\right\},\\
(-\infty,\beta]_{L^2(\Omega;H)}&:=&\left\{f \in L^2(\Omega;H): f(x)\le \beta\hbox{ for a.e. } x\in \Omega\right\},
\end{eqnarray*}
for $\alpha,\beta\in H$ with $-\alpha,\beta\in H_+$, so that these order intervals actually contain $0$. In particular, if $f\in H^1(\Omega;H)$, then $f^+,f^-\in H^1(\Omega;H)$.

2) Observe that although $f^+$ is formally given by the composition of a Lipschitz continuous mapping on $H$ and a function in $H^1(\Omega;H)$, providing a chain rule is not trivial as Rademacher's theorem fails to hold in infinite dimensional spaces and it is in particular not easy to understand in which sense the orthogonal projection of $H$ onto $H_+$ is ``differentiable a.e.'', as one would expect in the finite dimensional case. 

To this aim, let $f\in H^1(\Omega;H)$. By 1), one has in particular and by definition of $H^1(\Omega;H)$ that
$$\int_\Omega f^+(x)\nabla h(x) dx=-\int_\Omega \nabla f^+(x) h(x) dx\qquad\hbox{ for all }h\in C^\infty_c(\Omega;{\mathbb C})$$
in the sense of $H^n$-valued Bochner integrals. In other words, the above integrals define an element of $L^2(\Xi,{\mathbb C})^n$. Accordingly,
$$\left(\int_\Omega f^+(x)\nabla h(x) dx\right) (\xi)=-\left(\int_\Omega \nabla f^+(x) h(x) dx\right) (\xi)\qquad\hbox{ for all }h\in C^\infty_c(\Omega;{\mathbb C})\hbox{ and a.e. } \xi\in \Xi,$$
and therefore
$$\int_\Omega f^+(x,\xi)\nabla h(x) dx=-\int_\Omega \nabla f^+(x,\xi) h(x) dx\qquad\hbox{ for all }h\in C^\infty_c(\Omega;{\mathbb C})\hbox{ and a.e. } \xi\in \Xi:$$
this can be checked by first considering step functions and then going to the limit. We deduce that $f^+(\cdot,\xi)\in H^1(\Omega;{\mathbb C})$ for a.e.\ $\xi\in\Xi$. Since this is a scalar-valued function, we can apply the usual differentiation formula and deduce from~\cite[Lemma~7.6]{GilTru01} that
\begin{equation}\label{diffinal2}
\nabla f^+(\cdot,\xi)={\bf 1}_{\{f(\cdot,\xi)\ge 0\}} \nabla f(\cdot,\xi)\qquad \hbox{for a.e. } \xi\in \Xi.
\end{equation}
Now, because $f^+\in H^1(\Omega;H)$, the weak derivative of $f^+$ is necessarily given by~\eqref{diffinal2} outside a subset of $\Omega\times \Xi$ of zero measure.

We emphasize that the characteristic function is defined by means of subsets of $\Omega$ such that some inequality is satisfied by a \emph{scalar}-valued function. In fact the two subsets $\{f(\cdot,\xi)\ge 0\},\{f(\cdot,\xi)< 0\}$ define a partition of $\Omega$ for a.e.\ $\xi\in\Xi$.

3) We are now in the position to prove the main assertion. Since $u-v\in H^1(\Omega;H)$, we deduce from 2) that $(u-v)^+(\cdot,\xi),(u-v)^-(\cdot,\xi)\in H^1(\Omega;{\mathbb C})$ and the identities
\begin{equation*}\label{step1}
\nabla(u-v)^+(\cdot,\xi)=\mathbf 1_{\{u(\cdot,\xi)\ge v(\cdot,\xi)\}}(\nabla u-\nabla v)(\cdot,\xi),\qquad 
\nabla(u-v)^-(\cdot,\xi)=\mathbf 1_{\{u(\cdot,\xi)<v(\cdot,\xi)\}}(\nabla u-\nabla v)(\cdot,\xi)
\end{equation*}
hold for a.e.\ $\xi\in\Xi$. Accordingly, both 
$$P_{(-\infty,v]}u=\min\{u,v\}=u-(u-v)^+\qquad\hbox{and}\qquad P_{[v,+\infty)}u=\max\{u,v\}=v+(u-v)^+$$ 
belong to $H^1(\Omega;H)$ and~\eqref{gteine} follows. 

The remaining assertions are proven likewise.
\end{proof}

\begin{theo}\label{positive1}
Let $a\in H^1(\Omega\times \Xi;{\mathbb C})$ be such that ${\bf a}=(a,a_{|\Gamma})\in V_{\mathcal Y}$ and consider the unbounded order interval 
\begin{eqnarray*}
[{\bf a},+\infty)_{{\mathcal L}^2}&:=&[a,\infty)_{L^2(\Omega;H)}\times [a_{|\Gamma},\infty)_{L^2(\Gamma;H)}\\
&\cong &\left\{f \in L^2(\Omega\times\Xi;{\mathbb C}): a(x,\xi)\le f(x,\xi)\hbox{ for a.e. } (x,\xi)\in \Omega\times\Xi\right\}\\
&&\qquad\qquad\times \left\{g \in L^2(\Gamma\times\Xi;{\mathbb C}): a(z,\xi)\le g(z,\xi)\hbox{ for a.e. } (z,\xi)\in \Gamma\times\Xi\right\}
\end{eqnarray*}
Then $(e^{t\Delta_{{\mathcal Y},{\mathcal S}}})_{t\geq 0}$ leaves invariant $[{\bf a},+\infty)_{{\mathcal L}^2}$ if and only if 
\begin{enumerate}[(i)]
\item
$P_{\mathcal Y}[{\bf a},+\infty)_{{\mathcal L}^2}\subset [{\bf a},+\infty)_{{\mathcal L}^2}$ and additionally
\item the inequality
\begin{eqnarray*}
0&\ge &\int_\Xi \int_{\{a(\cdot,\xi)> f(\cdot,\xi)\}} \nabla a(x,\xi)\overline{(\nabla f-\nabla a)(x,\xi)}dx\; d\xi\\
&&\quad +\gamma\int_\Xi \int_{\{a_{|\Gamma}(\cdot,\xi)> f_{|\Gamma}(\cdot,\xi)\}} \nabla a(z,\xi) \overline{(\nabla f(z,\xi)-\nabla a(z,\xi))}d\sigma(z)\;d\xi +\left({\mathcal S} \max\{a_{|\Gamma},f_{|\Gamma}\}|(f_{|\Gamma}-a_{|\Gamma})^-\right)_{\mathcal Y}
\end{eqnarray*}
holds for all $f\in H^1(\Omega\times \Xi;{\mathbb R})$ such that ${\bf f}=(f,f_{|\Gamma})\in V_{\mathcal Y}$.
\end{enumerate}
\end{theo}

\begin{proof}
By Lemma~\ref{ouhalemma}, $(e^{t\Delta_{{\mathcal Y},{\mathcal S}}})_{t\geq 0}$ leaves invariant the order interval $[{\bf a},\infty)_{{\mathcal L}^2}$ if and only if the associated orthogonal projection $P_{[{\bf a},\infty)_{{\mathcal L}^2}}$ leaves invariant $V_{\mathcal Y}$ and moreover $a(P_{[{\bf a},+\infty)_{{\mathcal L}^2}}{\bf f},{\bf f}-P_{[{\bf a},+\infty)_{{\mathcal L}^2}}{\bf f})\ge 0$ for all $H_{\mathbb R}$-valued ${\bf f}\in V_{\mathcal Y}$. By Lemma~\ref{gteinelemma}, the first condition is satisfied if and only if $P_{[{\bf a},+\infty)_{\mathcal Y}}{\mathcal Y}\subset{\mathcal Y}$. By~\cite[Lemma~2.3]{ManVogVoi05} this is equivalent to $P_{\mathcal Y}[{\bf a},+\infty)_{\mathcal Y}\subset[{\bf a},+\infty)_{\mathcal Y}$.

The second criterion can be deduced applying Lemma~\ref{gteinelemma} and observing that for all ${\bf f}\in V_{\mathcal Y}$
\begin{eqnarray*}
a_{\mathcal S}(P_{[{\bf a},+\infty)_{{\mathcal L}^2}}{\bf f},{\bf f}-P_{[{\bf a},+\infty)_{{\mathcal L}^2}}{\bf f})&=&
-a_{\mathcal S}(\max\{{\bf a},{\bf f}\},({\bf f}-{\bf a})^-)\\
&=&-\int_\Omega \left( \mathbf 1_{\{a> f\}} \nabla a + \mathbf 1_{\{a\le f\}} \nabla f|\mathbf 1_{\{a> f\}}(\nabla f-\nabla a)\right)_{H^n}dx\\
&&\quad -\gamma\int_{\Gamma} \left( \mathbf 1_{\{a_{|\Gamma}> f_{|\Gamma}\}} \nabla_\Gamma a_{|\Gamma} |\mathbf 1_{\{a_{|\Gamma}> f_{|\Gamma}\}}(\nabla_\Gamma f_{|\Gamma}-\nabla_\Gamma a_{|\Gamma})\right)_{H^{n-1}}d\sigma(z)\\
&&\quad -\gamma\int_{\Gamma} \left( \mathbf 1_{\{a_{|\Gamma}\le f_{|\Gamma}\}} \nabla_\Gamma f_{|\Gamma} |\mathbf 1_{\{a_{|\Gamma}> f_{|\Gamma}\}}(\nabla_\Gamma f_{|\Gamma}-\nabla_\Gamma a_{|\Gamma})\right)_{H^{n-1}}d\sigma(z)\\
&&\quad -\left({\mathcal S} \max\{a_{|\Gamma},f_{|\Gamma}\}|(f_{|\Gamma}-a_{|\Gamma}^-)\right)_{\mathcal Y}\\
&=&-\int_\Omega\int_\Xi \left(\mathbf 1_{\{a(\cdot,\xi)> f(\cdot,\xi)\}} \nabla a(x,\xi) + \mathbf 1_{\{a(\cdot,\xi)\le f(\cdot,\xi)\}} \nabla f(x,\xi)\right)\\
&&\qquad\qquad\cdot\overline{\left(\mathbf 1_{\{a(\cdot,\xi)> f(\cdot,\xi)\}}(\nabla f-\nabla a)(x,\xi)\right)}d\xi \; dx\\
&&\quad -\gamma\int_{\Gamma}\int_\Xi \left( \mathbf 1_{\{a_{|\Gamma}(\cdot,\xi)> f_{|\Gamma}(\cdot,\xi)\}} \nabla a_{|\Gamma}(z,\xi)\right)\\
&&\qquad\qquad\cdot \overline{\left(\mathbf 1_{\{a_{|\Gamma}(\cdot,\xi)> f_{|\Gamma}(\cdot,\xi)\}}(\nabla_\Gamma f(z,\xi)-\nabla_\Gamma a(z,\xi)\right)} d\xi\; d\sigma(z)\\
&&\quad -\gamma\int_{\Gamma}\int_\Xi \left( \mathbf 1_{\{a_{|\Gamma}(\cdot,\xi)\le f_{|\Gamma}(\cdot,\xi)\}} \nabla_\Gamma f(z,\xi)\right)\\
&&\qquad\qquad\cdot \overline{\left(\mathbf 1_{\{a_{|\Gamma}(\cdot,\xi)> f_{|\Gamma}(\cdot,\xi)\}}(\nabla_{|\Gamma} f(z,\xi)-\nabla_{|\Gamma} a(z,\xi)\right)}d\xi d\sigma(z)\\
&&\quad -\left({\mathcal S} \max\{a_{|\Gamma},f_{|\Gamma}\}|(f_{|\Gamma}-a_{|\Gamma})^-\right)_{\mathcal Y}.
\end{eqnarray*}
Applying Fubini's theorem we obtain
\begin{eqnarray*}
a_{\mathcal S}(P_{[{\bf a},+\infty)_{{\mathcal L}^2}}{\bf f},{\bf f}-P_{[{\bf a},+\infty)_{{\mathcal L}^2}}{\bf f})
&=&-\int_\Xi \int_{\{a(\cdot,\xi)> f(\cdot,\xi)\}} \nabla a(x,\xi)\overline{(\nabla f-\nabla a)(x,\xi)}dx\; d\xi\\
&&\quad -\gamma\int_\Xi \int_{\{a_{|\Gamma}(\cdot,\xi)> f_{|\Gamma}(\cdot,\xi)\}} \nabla a(z,\xi) \overline{(\nabla f(z,\xi)-\nabla a(z,\xi))}d\sigma(z)\;d\xi \\
&&\quad -\left({\mathcal S} \max\{a_{|\Gamma},f_{|\Gamma}\}|(f_{|\Gamma}-a_{|\Gamma})^-\right)_{\mathcal Y}.
\end{eqnarray*}
This concludes the proof.
\end{proof}

Analogous assertions hold for the order intervals $(-\infty,{\bf b}]_{{\mathcal L}^2}$.

In general, condition (ii) in Theorem~\ref{positive1} will rarely be satisfied. An easy, yet relevant special case is clearly that of constant ${\bf a}$, i.e., ${\bf a}(x,\xi)\equiv \alpha$ for some $\alpha\in H$ and a.e.\ $(x,\xi)\in\Omega\times \Xi$. In this case, condition (ii) reduces to the condition
\begin{equation}\label{possemigr}
\left({\mathcal S} \max\{a_{|\Gamma},f_{|\Gamma}\}|(f_{|\Gamma}-a_{|\Gamma})^-\right)_{\mathcal Y}\le 0.
\end{equation}
Observe that if in addition ${\mathcal S}\in {\mathcal L}(L^2(\Gamma;H))$, then by Lemma~\ref{ouhalemma} the validity of condition (ii) in Theorem~\ref{positive1} is equivalent to the invariance of $[{\bf a},+\infty)_{L^2(\Gamma;H)}$ under the semigroup generated by ${\mathcal S}$. E.g., positivity of the semigroup corresponds to invariance of $[{\bf 0},\infty)_{{\mathcal L}^2}$, while ${\mathcal L}^\infty$-contractivity\footnote{ By this we mean contractivity with respect to the norm of $L^\infty(\Omega\times\Xi;{\mathbb C})\times L^\infty(\Gamma\times\Xi;{\mathbb C})$.} can be formulated in terms of simultaneous invariance of both order intervals $[-{\bf 1},\infty)_{{\mathcal L}^2},(-\infty,{\bf 1}]_{{\mathcal L}^2}$. 

For the sake of further reference we introduce the following locality assumptions.

\begin{assums}\label{locality}
There exist a closed subspace $Y$ of $H$, a closed convex subset $C_H$ of $H$ and an operator $S\in{\mathcal L}(H)$ such that
\begin{itemize}
\item ${\mathcal Y}=\{f\in L^2(\Gamma;H): f(z)\in Y \hbox{ for a.e. } z\in\Gamma\}$,
\item $C_{L^2(\Omega;H)}:=\{f\in L^2(\Omega;H): f(x)\in C_H \hbox{ for a.e. } x\in\Omega\}$,
\item $C_{\mathcal Y}:=\{f\in {\mathcal Y}: f(z)\in C_H \hbox{ for a.e. } z\in\Gamma\}$,
\item $C_{{\mathcal L}^2}:=C_{L^2(\Omega;H)}\times C_{\mathcal Y}$ and 
\item ${\mathcal S}g= S\circ g$ for all $g\in H^\frac{1}{2}(\Gamma;H)$.
\end{itemize}
Moreover, $0\in C$ or else both $\Omega$ and $\Gamma$ have finite measure.
\end{assums}

Observe that under Assumptions~\ref{locality} the abstract problem $\rm(AV)$ becomes a parabolic problem with dynamic boundary condition
\begin{equation*}\tag{DBC}
\left\{\begin{array}{rcll}
\frac{\partial}{\partial t}{u}(t,x)&=& \Delta u(t,x),\qquad &t\geq 0,\; x\in \Omega,\\
u(t,z)&\in& Y, &t\geq 0,\; z\in\Gamma,\\
\frac{\partial }{\partial t}u(t,z)&=&P_Y\left(-\frac{\partial u}{\partial\nu}(t,z)+\left(\gamma \Delta_\Gamma-S\right) u(t,z)\right), &t\geq 0,\; z\in\Gamma,\\
u(0,x)&=&u_0(x),&x\in\Omega,\\
u(0,z)&=&v_0(z),&z\in\Gamma.
\end{array}
\right.
\end{equation*}
If $\Omega=(0,\infty)\times {\mathbb R}^{n-1}$, it is common in the literature to refer to this problem as ``diffusion on an open book'' (with dynamic boundary conditions). If $n=1$, this is nothing but the semi-infinite star considered in Example~\ref{firstexa}

Under the Assumptions~\ref{locality}, $\mathcal Y$ is a closed subspace of $L^2(\Gamma;H)$ and $C_{L^2(\Omega;H)},C_{\mathcal Y},C_{{\mathcal L}^2}$ are closed and convex subsets of $L^2(\Omega;H)$, $\mathcal Y$, and ${\mathcal L}^2$, respectively. With an abuse of notation we then write ${\mathcal Y}\equiv Y$, ${\mathcal S}\equiv S$ and $\Delta_{Y,S}$ instead of $\Delta_{{\mathcal Y},{\mathcal S}}$.

It is crucial that whenever Assumptions~\ref{locality} hold the orthogonal projections of $L^2(\Omega;H)$ onto $P_{L^2(\Omega;H)}$, of $\mathcal Y$ onto $P_{\mathcal Y}$ and hence of ${\mathcal L}^2$ onto $C_{{\mathcal L}^2}$ satisfy
\begin{eqnarray*}
P_{C_{L^2(\Omega;H)}} f &=& P_{C_H}\circ f\qquad \hbox{for all }{f}\in L^2(\Omega;H),\\
P_{C_{\mathcal Y}} g &=& P_{C_H}\circ g\qquad \hbox{for all }{g}\in {\mathcal Y},\\
P_{C_{{\mathcal L}^2}} {\bf f} &=& \begin{pmatrix} P_{C_H}\circ f\\ P_{C_H}\circ g\end{pmatrix} 
\qquad \hbox{for all }{\bf f}:=\begin{pmatrix}f\\ g\end{pmatrix}\in{\mathcal L}^2.
\end{eqnarray*}
The fact that the projections onto the above subsets of vector-valued function spaces are the compositions of a Lipschitz continuous mapping (namely, the projection $P_{C_H}$) and a function of class $H^1$ permits to apply the version of a chain rule obtained in Lemma~\ref{gteinelemma}. Furthermore, due to the local structure of the sets $C_{L^2(\Omega;H)}$ and $C_{\mathcal Y}$, one sees that in particular
$$(P_{C_{L^2(\Omega;H)}}\circ f)_{|\Gamma}= (P_{C_{\mathcal Y}}\circ f_{|\Gamma})\qquad\hbox{for all }{\bf f}\in V_{\mathcal Y}.$$

\begin{theo}\label{invsubsp}
Impose the Assumptions~\ref{locality}. Then $C_{{\mathcal L}^2}$ is left invariant under $(e^{t\Delta_{{\mathcal Y},{\mathcal S}}})_{t\ge 0}$ if and only if
\begin{enumerate}[(i)]
\item the inclusion $P_Y C_H \subset C_H$ holds and additionally
\item the semigroup $(e^{-tS})_{t\ge 0}$ leaves $C_H$ invariant.
\end{enumerate}
\end{theo}

Comparable results have been obtained in the context of networks in~\cite{KosPotSch08,KanKlaVoi09}.

\begin{proof}
First of all, we show that the inclusion $P_{C_{{\mathcal L}^2}} V_{\mathcal Y} \subset V_{\mathcal Y}$ holds if and only if the inclusion $P_{\mathcal Y} C_{\mathcal Y} \subset C_{\mathcal Y}$ holds. Orthogonal projections onto closed convex subsets of a Hilbert space are Lipschitz continuous mappings, hence as already observed by~\cite[Lemma~7.3]{CarMug09} $P_{C_{{\mathcal L}^2}}$ maps $H^1(\Omega;H)\times H^1(\Gamma;H)$ into itself -- i.e., the weak differentiability conditions is satisfied independently of the boundary conditions.
Consequently, $P_{C_{{\mathcal L}^2}}V_{\mathcal Y}\subset V_{\mathcal Y}$ if and only if $f_{|\Gamma}\in {\mathcal Y}$ implies $P_{C_{\mathcal Y}} f_{|\Gamma}\in \mathcal Y$, for all ${\bf f}\in H^1(\Omega;H)$. The proof can be completed reasoning as in~\cite[Prop.~4.2]{CarMug09}.

By Lemma~\ref{ouhalemma}, invariance of $C_{\mathcal Y}$ under $(e^{t\Delta_{{\mathcal Y},{\mathcal S}}})_{t\ge 0}$ is now equivalent to $P_{\mathcal Y} C_{\mathcal Y} \subset C_{\mathcal Y}$ and 
$${\rm Re}a_{\mathcal S}(P_{C_{{\mathcal L}^2}}{\bf f}, (I-P_{C_{{\mathcal L}^2}}){\bf f})\ge 0\qquad\hbox{for all }{\bf f}\in V_{\mathcal Y}.$$ 
Due to locality of the forms associated with the Laplacian on $\Omega$ and the Laplace--Beltrami operator on $\Gamma$ (and hence of the form $(a_{\mathcal S},V_{\mathcal Y})$), a direct computation shows that
$${\rm Re}a_{\mathcal S}(P_{C_{{\mathcal L}^2}}{\bf f}, (I-P_{C_{{\mathcal L}^2}}){\bf f})={\rm Re}(SP_{C_\mathcal Y} f_{|\Gamma}|(I-P_{C_\mathcal Y})f_{|\Gamma})_{\mathcal Y}.$$
By density, the latter term is $\ge 0$ for all ${\bf f}\in V_{\mathcal Y}$ if and only if
$${\rm Re}(SP_{C_{\mathcal Y}}g|(I-P_{C_{\mathcal Y}})g)_{\mathcal Y}\ge 0\qquad \hbox{for all }g\in {\mathcal Y}.$$
By a localisation argument this is equivalent to asking that
$${\rm Re}(SP_{C_H}x|(I-P_{C_H}))x)_H\ge 0\qquad \hbox{for all }x\in H.$$
A further application of Lemma~\ref{ouhalemma} concludes the proof, since $(S\cdot|\cdot)_H$ is the form associated with $-S$.
\end{proof}

In the previous theorem, it is not too restrictive to consider sets of the form $C_{L^2(\Omega;H)}\times C_{L^2(\Gamma;H)}$ -- i.e., to restrict ourselves to study invariance of sets of those functions pointwise belonging to the same subset of $H$, both on $\Omega$ and on the boundary $\Gamma$. In fact, the following holds.

\begin{prop}\label{C=D}
Let $C,D\subset H$ be closed convex subsets. If $C_{L^2(\Omega;H)}\times D_{L^2(\Gamma;H)}$ is invariant under $(e^{t\Delta_{{\mathcal Y},{\mathcal S}}})_{t\ge 0}$, then $C=D$.
\end{prop}

\begin{proof}
We only consider the case of $\Omega,\Gamma$ with bounded measure. The general case will then follow by localisation arguments. Let first $C\not\subset D$, say $v\in C\setminus D$. Take ${\bf f}\in V_{\mathcal Y}$ such that $f=1_\Omega\otimes v$ -- i.e., $f\equiv v$: then $f\in C_{L^2(\Omega;H)}$ and $f_{|\Gamma}=1_{\Gamma} \otimes v\not\in D_{L^2(\Gamma;H)}$. Then $$P_{C_{L^2(\Omega;H)}\times D_{L^2(\Gamma;H)}}{\bf f}=\begin{pmatrix} 1\otimes v\\ 1\otimes P_{D} v\end{pmatrix},$$
i.e., $(1_\Omega\otimes v)_{|\Gamma}\not= 1_{\Gamma}\otimes P_{D} v$ and accordingly 
$P_{C_{L^2(\Omega;H)}\times D_{L^2(\Gamma;H)}}{\bf f}\not\in V_{\mathcal Y}$. The case of $D\not\subset C$ can be treated likewise.
\end{proof}

We mention that domination of semigroups can also be discussed. E.g., the following can be shown mimicking the proof of~\cite[Cor.~2.22]{Ouh05}. This results extends~\cite[Prop.~2.8]{AreMetPal03} and~\cite[Prop.~4.2]{MugRom07}.

\begin{prop}\label{dominvector}
Impose the Assumptions~\ref{locality} and let $P_{\mathcal Y}$ be a positive operator. Let $S_1, S_2$ be $L^\infty(\Gamma;{\mathcal L}_s(H))$-functions\footnote{ Here, we denote by $L^\infty(\Gamma;{\mathcal L}_s(H))$ the space of all measurable and essentially bounded functions from $\Gamma$ to ${\mathcal L}(H)$ with respect to the strong operator topology.}. Define operators ${\mathcal S}_1,{\mathcal S}_2$ by
$${\mathcal S}_1g= S_1\circ g\quad\hbox{ and }\quad{\mathcal S}_2g= S_2\circ g\qquad\hbox{ for all } g\in H^\frac{1}{2}(\Gamma;H).$$
 Consider two sesquilinear forms $a_1,a_2$ defined by
\begin{equation*}
a_1 \left({\bf f},{\bf g}\right):=\int_\Omega \left(\nabla f(x) | \nabla g(x)\right)_{H^n} dx+\gamma\int_{\Gamma} \left (\nabla f(z)|\nabla g(z)\right)_{H^{n-1}} d\sigma(z)+( S_1 f_{|\Gamma}| g_{|\Gamma})_{\mathcal Y}
\end{equation*}
and
\begin{equation*}
a_2 \left({\bf f},{\bf g}\right):=\int_\Omega \left(\nabla f(x) | \nabla g(x)\right)_{H^n} dx+\gamma\int_{\Gamma} \left (\nabla f(z)|\nabla g(z)\right)_{H^{n-1}} d\sigma(z)+( S_2 f_{|\Gamma}| g_{|\Gamma})_{\mathcal Y},
\end{equation*}
both defined on $V_Y$, and the associated operators $\Delta_{Y,S_1},\Delta_{Y,S_2}$. Then the following assertions hold.
\begin{enumerate}
\item The semigroup $(e^{t\Delta_{Y,S_1}})_{t\geq 0}$ is dominated by $(e^{t\Delta_{Y,S_2}})_{t\geq 0}$, i.e.
$$|e^{t\Delta_{Y,S_1}}f(x,\xi)| \le e^{t\Delta_{Y,S_2}}|f|(x,\xi),\qquad t\ge 0,\; f\in L^2(\Omega\times \Xi;{\mathbb C}),\; x\in\Omega,\; \xi\in\Xi,$$
 if and only if 
$${\rm Re}\;( S_1 f_{|\Gamma}| g_{|\Gamma})_{\mathcal Y} \ge ( S_2 |f|_{|\Gamma}| |g|_{|\Gamma})_{\mathcal Y}$$
 for all $u,v\in V_Y$ such that $u\overline{v}\ge 0$.

\item Let $ S_1(z), S_2(z)$ be positive operators for a.e.\ $z\in \Gamma$. Then the semigroup $(e^{t\Delta_{Y,S_1}})_{t\geq 0}$ is dominated by $(e^{t\Delta_{Y,S_2}})_{t\geq 0}$
%, i.e.
%$$|e^{t\Delta_{Y,S_1}}f(x,\xi)| \le e^{t\Delta_{Y,S_2}}|f|(x,\xi),\qquad t\ge 0,\; f\in L^2(\Omega\times \Xi;{\mathbb C}),\; x\in\Omega,\; \xi\in\Xi,$$
 if and only if $S_1(z)-  S_2(z)$ is a positive operator for a.e.\ $z\in\Gamma$.
\end{enumerate}

\end{prop}

\begin{rem}
In the usual theory of semigroup domination, both the dominating and the dominated semigroup have to act on the same space, or else one of them has to act on a space of scalar-valued functions, see~\cite{ManVogVoi05} and references therein. This rules out several interesting case in our context, due to the fact the boundary conditions also determine the state space -- and hence semigroups governing equations with different boundary conditions cannot been compared. E.g., it would be natural to expect that all semigroups $(e^{t\Delta_{{\mathcal Y},{\mathcal S}}})_{t\ge 0}$ dominate the semigroup that governs the heat equation with (uncoupled) Dirichlet boundary conditions, provided that condition~\eqref{possemigr} holds.
\end{rem}

While it is known that many relevant properties are shared by the heat equation with either non-dynamic or dynamic boundary conditions, to the best of our knowledge a structural relation between these phenomena had not yet been observed. The following is a straightforward consequence of Theorem~\ref{invsubsp} and~\cite[Prop.~4.3]{CarMug09}.

\begin{cor}\label{almosthere2}
Impose the Assumptions~\ref{locality}. Then $C_{{\mathcal L}^2}$ is left invariant under $(e^{t\Delta_{Y,S}})_{t\ge 0}$ if and only if
$C_{L^2(\Omega;H)}$ is left invariant under the semigroup governing the parabolic problem
\begin{equation}\tag{NDBC}
\left\{\begin{array}{rcll}
\frac{\partial}{\partial t}{u}(t)&=& \Delta u(t),\qquad &t\geq 0,\\
u(t)_{|\Gamma}&\in& Y, &t\geq 0,\\
\frac{\partial u(t)}{\partial\nu}+Su(t)_{|\Gamma}&\in& Y^\perp, &t\geq 0.\\
u(0)&=&u_0,
\end{array}
\right.
\end{equation}
with time-independent boundary conditions.
\end{cor}

Observe that the semigroup governing $\rm(NDBC)$ is generated by the operator associated with $a_{\mathcal S}$ (with $\gamma=0$), but considered as a sesquilinear form acting on the Hilbert space $\{f\in H^1(\Omega;H):f_{|\Gamma}\in \mathcal Y\}\hookrightarrow L^2(\Omega;H)$ rather than $V_{\mathcal Y}\hookrightarrow {\mathcal L}^2$, cf.~\cite{CarMug09}.

\begin{exa}
As shown in~\cite{FavGolGol02,AreMetPal03}, remarkable properties of the (scalar-valued) heat equation with Wentzell--Robin (dynamic) boundary conditions include positivity and contractivity with respect to the $\infty$-norm of the semigroup that governs it. In the light of Corollary~\ref{almosthere2}, these properties actually follow from the same properties enjoyed by the heat equation with corresponding Robin (time-independent) boundary conditions.
\end{exa}

Observe in particular that
$${\mathcal L}^p\equiv{\mathcal L}^p_{\mathcal Y}:=L^p(\Omega;H)\times (L^p(\Gamma;H)\cap{\mathcal Y}),\qquad p\in [1,\infty],$$
are Bochner spaces with respect to a suitable product measure. Assume both $(e^{t\Delta_{{\mathcal Y},{\mathcal S}}})_{t\ge 0}$ and its adjoint to be ${\mathcal L}^\infty$-contractive: under the Assumptions~\ref{locality} this can be characterized by means of Theorem~\ref{invsubsp}, with $C_H=(-\infty,{\bf 1}]_H\cap [{\bf 1},\infty)_H$.

\begin{cor}
Assume both $(e^{t\Delta_{{\mathcal Y},{\mathcal S}}})_{t\ge 0}$ and its adjoint to be ${\mathcal L}^\infty$-contractive and let $n\ge 2$. Then $(e^{t\Delta_{{\mathcal Y},{\mathcal S}}})_{t\ge 0}$ extrapolates to a consistent family of operator semigroups on ${\mathcal L}^p$, $p\in [1, \infty]$. These semigroups are strongly continuous and analytic for $p\in(1,\infty)$. 

Moreover, $(e^{t\Delta_{{\mathcal Y},{\mathcal S}}})_{t\ge 0}$ is ultracontractive, i.e., it satisfies the estimate
\begin{equation*}
\| e^{t\Delta_{{\mathcal Y},{\mathcal S}}}{\bf f}\|_{\mathcal{L}^\infty} \leq
M_\mu t^{-\frac{\mu}{2}}\|{\bf f}\|_{\mathcal{L}^2}
\quad\hbox{ for all }t\in (0,1],\; {\bf f}\in{\mathcal{L}^2}
\end{equation*}
where
\begin{equation*}
\mu\in
\begin{cases}
[n-1,\infty), &\hbox{ if }\; n\geq 3,\\
(1,\infty), &\hbox{ if }\; n=2,
\end{cases}
\end{equation*}
and some constant $M_\mu$. The same estimates are satisfied by the dual semigroup.
\end{cor}

Additional conditions ensuring strong continuity for $p=1$ are known, cf.~\cite[{\S}7.2.1]{Are06} for the scalar case.

\begin{proof}
The assertion on extrapolation follows applying a vector-valued version of Riesz--Thorin's interpolation theorem, cf.~\cite[p.~77]{JohLin01}.
The second assertion can be proved as in the scalar-valued case, applying a known characterisation of ultracontractivity (see~\cite[{\S}12.2]{Are06}) based on standard Sobolev embeddings, cf.~\cite[Lemma~3.8]{MugRom06}. It can be easily seen that all the involved techniques carry over to the vector-valued case.
\end{proof}

\begin{rems}
1) By~\cite[Lemma~3.3]{MugNit09}, $(e^{t\Delta_{{\mathcal Y},{\mathcal S}}})_{t\ge 0}$ consists of kernel operators for all $t>0$.

2) It is remarkable that the above mentioned criterion for ultracontractivity based on Sobolev embeddings only applies if $n>1$. In the scalar case, a common workaround is to deduce ultracontractivity from the Nash inequality. Unfortunately, the Nash inequality seems to extend to the vector-valued case only if the space $H$ is finite dimensional. This is why we are not able to prove the above result in the case of $n=1$ -- which in particular corresponds to the relevant case of networks with infinitely many edges.

% 3) What about generation in a space of continuous functions? If $(e^{t\Delta_{{\mathcal Y},{\mathcal S}}})_{t\ge 0}$ is submarkovian and $\Omega$ has finite measure and $C^\infty$-boundary, then the semigroup maps ${\mathcal L}^2$ into
% $$C(\overline{\Omega};H)\equiv\{{\bf f}\in C(\overline{\Omega};H)\times C(\Gamma;H):f_{|\Gamma}\in{\mathcal Y}\},$$ 
% then one can mimick the proof of~\cite[Prop.~3.2]{AreMetPal03} and deduce that the part of $\Delta_{{\mathcal Y},{\mathcal S}}$ in $C(\overline{\Omega};H)\times C(\Gamma;H)$ is a resolvent positive operator, hence by~\cite[Thm.~3.11.9]{AreBatHie01} it generates a positive strongly-continuous semigroup on the closure of its domain provided that $\Omega,\Gamma$ have finite measure. However, we have seen that submarkovianity of the semigroup is a very strong assumption in our vector-valued setting.
\end{rems}

A semigroup on an $L^2$-space is said to be irreducible if the only closed ideals of $L^2$ left invariant under the semigroup are the trivial ones. If $Y$ is a closed ideal of $H$, then clearly $(e^{t\Delta_{Y,0}})_{t\ge 0}$ leaves invariant $L^2(\Omega;Y)\times L^2(\Gamma;Y)$, which is a closed ideal of ${\mathcal L}^2$. Thus, uncoupled boundary conditions jeopardize irreducibility. 

More generally, we observe that if $\mathcal P:\Omega\to {\mathcal L}(H)$ is a strongly measurable function such that ${\mathcal P}(x)$ is an orthogonal projection onto a closed ideal of $H$ for a.e.\ $x\in \Omega$, then the subspace
\begin{equation}
I_\mathcal P := \{f \in L^2(\Omega;H): f(x) \in {\rm Range \,}\mathcal P(x) \hbox{ for a.e. }  x \in \Omega \}
\end{equation}
is a closed ideal of $L^2(\Omega;H)$, too. In fact,  all closed ideals of $L^2(\Omega;H)$ are of this form, as it is proven in~\cite{CarMug09b}. Similarly, if the Assumptions~\ref{locality} hold one can see that each closed ideals of ${\mathcal L}^2$ is the range of an operator-valued strongly measurable mapping $\mathcal P$ defined on the product measure space $\Omega\oplus \Gamma$ and such that
\begin{itemize}
 \item ${\mathcal P}(x)$ is an orthogonal projection onto a closed ideal of $H$ for a.e.\ $x\in \Omega$ and
 \item ${\mathcal P}(z)$ is an orthogonal projection onto a closed ideal of $Y$ for a.e.\ $z\in\Gamma$.
\end{itemize}

\begin{prop}
Impose the Assumptions~\ref{locality}. Then $(e^{t\Delta_{Y,0}})_{t\geq 0}$ is irreducible if and only if $P_Y$ is irreducible and $\Omega$ is connected.
\end{prop}

Observe that in the scalar case $H={\mathbb C}$ the orthogonal projections on both subspaces of $H$ are irreducible.

\begin{proof}
It is clear that the semigroup is not irreducible if $\Omega$ is unconnected, since it lets invariant the closed ideals consisting of those functions supported in any of the connected components.

Let now $P_Y$ be non-irreducible, i.e., consider a non-trivial closed ideal $J_H$ of $H$ such that $P_Y J_H \subset J_H$. Then by Theorem~\ref{invsubsp} we conclude that $J_{L^2(\Omega;H)}\times J_{\mathcal Y}$ is a closed ideal of ${\mathcal L}^2$ that is left invariant under the semigroup, i.e., $(e^{t\Delta_{Y,0}})_{t\geq 0}$ is not irreducible.

Let conversely $(e^{t\Delta_{Y,0}})_{t\geq 0}$ be non-irreducible. Then there exists a non-trivial closed ideal of ${\mathcal L}^2$ that is invariant under $(e^{t\Delta_{Y,0}})_{t\geq 0}$. By Proposition~\ref{C=D} such an ideal is necessarily of the form $C_{L^2(\Omega;H)}\times C_{L^2(\Gamma;H)}$. Now, we can apply Theorem~\ref{invsubsp} and deduce the claim.
\end{proof}

\begin{rem}
In the scalar case, it is known that irreducibility is equivalent to a strong parabolic maximum principle, provided that the semigroup is positive, cf.~\cite[{\S}2.2]{Ouh05} -- but this characterisation fails to hold in the general vector-valued case. E.g., the heat semigroup $(e^{t\Delta})_{t\ge 0}$ on $L^2({\mathbb R};{\mathbb R}^2)$ is not irreducible because $L^2({\mathbb R};{\mathbb R}\times\{0\})$ is a non-trivial closed ideal left invariant under the semigroup. However, it does map nonzero positive functions $f$ to functions $e^{t\Delta}f$ satisfying $e^{t\Delta}f(x)>0$\footnote{ I.e., $e^{t\Delta}f(x)$ is a nonzero, positive vector of ${\mathbb R}^2$.} for all $t>0$ and a.e.\ $x\in\mathbb R$. 
\end{rem}

\section{An example: Diffusion on a star-shaped network}\label{examplesect}

Throughout this section we consider the setting presented in Example~\ref{firstexa}. Observe that the Assumptions~\ref{locality} are satisfied whenever we discuss invariance of a set $C_H$ which is either a subspace or an order interval containg $0$. We are going to present some interesting behaviour even in this elementary setting. Actually, same properties hold for more general diffusion on domains, rather than intervals. Also, by Corollary~\ref{almosthere2} all the results in this section hold for the semigroups governing $\rm(NDBC)$ and $\rm(DBC)$ alike. Thus, we explicitly refer to the case of time-independent boundary conditions only.

It has been proved in~\cite[{\S}5]{CarMug09} that the semigroup governing $\rm(NDBC)$ is positive if $Y=\langle {\mathbf 1}\rangle$ (i.e., under so-called Kirchhoff boundary conditions) and not positive if $Y=\langle {\mathbf 1}\rangle^\perp$ (i.e., under so-called \emph{anti-Kirchhoff boundary conditions} as considered e.g.\ in~\cite{Kuc04,FulKucWil07,Pos07,AlbCacFic07}), provided that $-S$ generates a positive semigroup on $Y$ (i.e., $-S$ is a real matrix with positive off-diagonal entries).

Similarly, assume that $-S$ generates an $L^\infty$-contractive semigroup on $Y$ and that $H={\mathbb C}^N$. Then by~\cite[Lemma~6.1]{Mug07}, this can be characterized by the fact that the entries $s_{ij}$ of $S$ satisfy
$$\sum_{j\not=i}|s_{ij}|\le {\rm Re}s_{ii}\qquad \hbox{for all }i,$$
cf.\ also~\cite[Rem.~3.8.(2)]{CarMug09}. Then one can prove that the heat semigroup is $L^\infty$-contractive under Kirchhoff boundary conditions for all $N\in\mathbb N$, whereas in the anti-Kirchhoff case it is $L^\infty$-contractive if and only if $N=2$.

\bigskip
For the sake of simplicity, in the remainder of this section we let $S=0$. 

A semi-infinite star with two edges can be identified with a line. More precisely, up to the canonical isometric isomorphism $U$ defined by
$$(Uf)(x):=\begin{pmatrix}
f(x)\\ f(-x)
  \end{pmatrix},\qquad x\ge 0,$$
functions in $L^2({\mathbb R};{\mathbb C})$ and in $L^2((0,+\infty);{\mathbb C}^2)$ may be identified. Accordingly, a function $(f_1,f_2)\in L^2((0,+\infty);{\mathbb C}^2)$ is called \emph{even} (resp., \emph{odd}) if $f_1(x)=f_2(x)$ (resp., if $f_1(x)+f_2(x)=0$) for a.e.\ $x\in (0,+\infty)$. More generally, we call a function $f\in L^2(\Omega;\mathbb C^N)$ \emph{even} (resp., \emph{odd}) if $f(x)\in \langle{\mathbf 1}\rangle$ (resp., if $f(x)\in \langle{\mathbf 1}^\perp\rangle$) for a.e.\ $x\in\Omega$. 

By Theorem~\ref{invsubsp} both the diffusion semigroups with Kirchhoff (i.e., $Y=\langle{\bf 1}\rangle$) and anti-Kirchhoff (i.e., $Y=\langle{\bf 1}\rangle^\perp$) boundary conditions leave invariant the set of even functions as well as the set of odd ones. If $N=2$, then it is easy to see that these are in fact the \emph{only} boundary conditions leading to invariance of any of these both sets.

\bigskip
Now, consider a semi-infinite star with only two edges, i.e., $H={\mathbb C}^2$. Neglecting the trivial (uncoupled) boundary conditions defined by $Y=\{0\}$ and $Y={\mathbb C}^2$ we can consider all $1$-dimensional subspaces ${\mathcal Y}\equiv Y_\xi$ of ${\mathbb C}^2$ by means of the parametrisation
\begin{equation*}\label{parmatrix}
P_{Y_\xi}:=\begin{pmatrix}
\cos^2 \xi & \sin \xi\; \cos \xi\\
\sin \xi\; \cos \xi & \sin^2 \xi
\end{pmatrix},\qquad \xi\in [0,\pi),
\end{equation*}
where ${Y_\xi}$ denotes the range of the orthogonal projection $P_{Y_\xi}$. Observe that $\xi=0$, $\xi=\frac{\pi}{4}$, $\xi=\frac{\pi}{2}$ and $\xi=\frac{3\pi}{4}$ correspond to uncoupled Dirichlet/Neumann, to Kirchhoff, to uncoupled Neumann/Dirichlet and to anti-Kirchhoff boundary conditions, respectively, as can be checked directly.

We are going to discuss the submarkovian property of the semigroup associated with these subspaces in dependence of $\xi$. A direct computation shows that the semigroup $(e^{t\Delta_{Y_\xi,0}})_{t\ge 0}$ is positive if and only if $\xi \in [0,\frac{\pi}{2}]$. Furthermore, by Theorem~3.6 the semigroup that governs $\rm(NDBC)$ is $L^\infty(\Omega\times \Xi;{\mathbb C})$-contractive if and only if $P_{Y_\xi}$ is $L^\infty(\Xi;{\mathbb C})$-contractive
%\footnote{ I.e., contractive on the Hilbert lattice $H$ endowed with the $\infty$-norm.}
, i.e., if and only if the inequalities
$$\cos^2 \xi+|\sin\xi\; \cos\xi|\le 1\qquad\hbox{and}\qquad |\sin\xi\; \cos\xi|+\sin^2 \xi\le 1$$
hold simultaneously. The former (resp., the latter) inequality holds if and only if $\xi\not\in (0,\frac{\pi}{4})\cup(\frac{3\pi}{4},\pi)$ (resp., if and only if $\xi\not\in (\frac{\pi}{4},\frac{\pi}{2}),(\frac{\pi}{2},\frac{3\pi}{4})$).

\begin{center}
%% \begin{tikzpicture}[domain=0:2*pi]
%% \draw[->, semithick] (0,0) -- (6.8,0) node[below left] {$\xi$};
%% \draw[->, semithick]  (0,-1.2) -- (0,1.5) node[above] {};
%% \draw (0,-1.1) -- (0,1.5);
%% \draw (pi/2,-1.1) -- (pi/2,1.5);
%% \draw (pi,-1.1) -- (pi,1.5);
%% \draw (3*pi/2,-1.1) -- (3*pi/2,1.5);
%% \draw (2*pi,-1.1) -- (2*pi,1.5);
%% \draw (0,1) -- (2*pi,1);
%% \draw (0,-1) -- (2*pi,-1);
%% \draw (pi/2,-.1) -- (pi/2,.1) node[below left=1pt] {$\frac{\pi}{4}$};
%% \draw (pi,-.1) -- (pi,.1) node[below left=1pt] {$\frac{\pi}{2}$};
%% \draw (3*pi/2,-.1) -- (3*pi/2,.1) node[below left=1pt] {$\frac{3\pi}{4}$};
%% \draw (2*pi,-.1) -- (2*pi,.1) node[below left=1pt] {$\pi$};
%% \foreach \y in {-1,0,1} \draw (-.1,\y) -- (.1,\y) node[left] {$\y$};
%% \draw[color=violet] plot function{cos(.5*x)*cos(.5*x)+abs(sin(.5*x)*cos(.5*x))} node[right=12pt] {$\cos^2(\xi)+|sin(\xi)\;\cos(\xi)|$};
%% \draw[color=blue] plot function{sin(.5*x)*sin(.5*x)+abs(sin(.5*x)*cos(.5*x))} node[right=12pt] {$\sin^2(\xi)+|sin(\xi)\;\cos(\xi)|$};
%% \draw[color=gray] plot function{sin(.5*x)*cos(.5*x)} node[below right=12pt] {$\sin \xi \; \cos \xi$};
%% \end{tikzpicture}\\
\begin{figure}[htbp]
    \scalebox{.8}
{\includegraphics{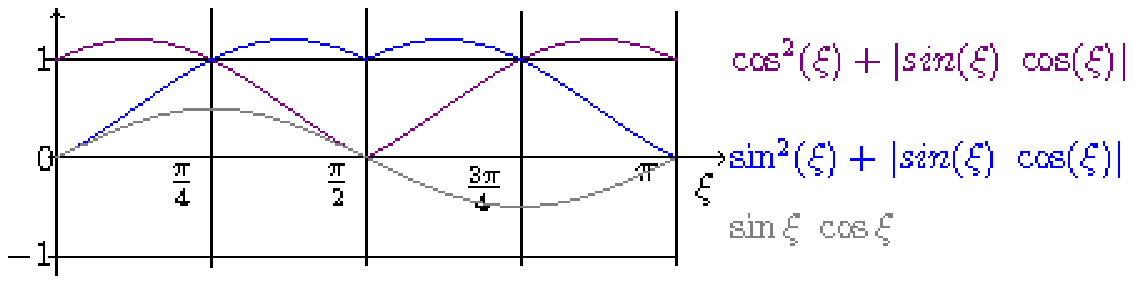}}
\end{figure}
{\sc Figure 1}
\end{center}

Therefore, the $L^\infty$-contractivity of the semigroup associated with Kirchhoff boundary conditions  represents a singularity. In particular, a submarkovian semigroup is generated \emph{exactly} in the following five cases:
\begin{itemize}
 \item with uncoupled Dirichlet/Dirichlet boundary conditions,
\item with uncoupled Neumann/Neumann boundary conditions,
\item with uncoupled Dirichlet/Neumann boundary conditions,
\item with uncoupled Neumann/Dirichlet boundary conditions and finally
\item with Kirchhoff boundary conditions.
\end{itemize}
Similarly, we can consider general boundary conditions defined by $1$-dimensional subspaces of $H$ for a semi-infinite star with $3$ edges ($H={\mathbb C}^3$). They can be investigated by means of spherical boundary conditions, i.e., considering spaces ${\mathcal Y}\equiv Y_{\xi,\phi}$ that are ranges of the orthogonal projections
\begin{equation}\label{bcmatrix}
P_Y\equiv P_{Y_{\xi,\phi}}=\begin{pmatrix}
\sin^2 \xi\; \cos^2\phi & \sin^2 \xi\; \sin\phi\; \cos\phi& \sin\xi \; \cos\xi\; \cos\phi\\
\sin^2 \xi\; \sin\phi\; \cos\phi & \sin^2\xi \;\sin^2\phi & \sin\xi\; \cos\xi\;\sin\phi\\
\sin\xi \; \cos\xi\; \cos\phi & \sin\xi\; \cos\xi\;\sin\phi & \cos^2 \xi\end{pmatrix},\qquad \xi,\phi\in [0,2\pi). 
\end{equation}
Analysing the behaviour of $P_{Y_{\xi,\phi}}$ in dependence of $\xi,\phi$ as done above for $P_{Y_\xi}$ is less elementary. While the matrix is clearly positive if and only if $\xi,\phi\in [0,\frac{\pi}{2}]\cup[\pi,\frac{3\pi}{2}]$, it is not clear how to determine all the values $\xi,\phi$ leading to $L^\infty$-contractivity, i.e., all the values $\xi,\phi$ such that the three functions
\begin{eqnarray*}
&&\sin^2 \xi\; \cos^2\phi +| \sin^2 \xi\; \sin\phi\; \cos\phi|+| \sin\xi \; \cos\xi\; \cos\phi|,\\
&&|\sin^2 \xi\; \sin\phi\; \cos\phi| + \sin^2\xi \;\sin^2\phi +| \sin\xi\; \cos\xi\;\sin\phi|\quad \hbox{ and}\qquad\qquad \xi,\phi\in [0,\pi)\\
&&|\sin\xi \; \cos\xi\; \cos\phi|+ |\sin\xi\; \cos\xi\;\sin\phi|+ \cos^2 \xi
\end{eqnarray*}
are simultaneously $\le 1$, corresponding to the three conditions for $L^\infty$-contractivity associated with the three rows of the matrix $P_{Y_{\xi,\phi}}$ in~\eqref{bcmatrix}.
\begin{center}
\begin{figure}[htbp]
  \centering
    \scalebox{.6}{\includegraphics{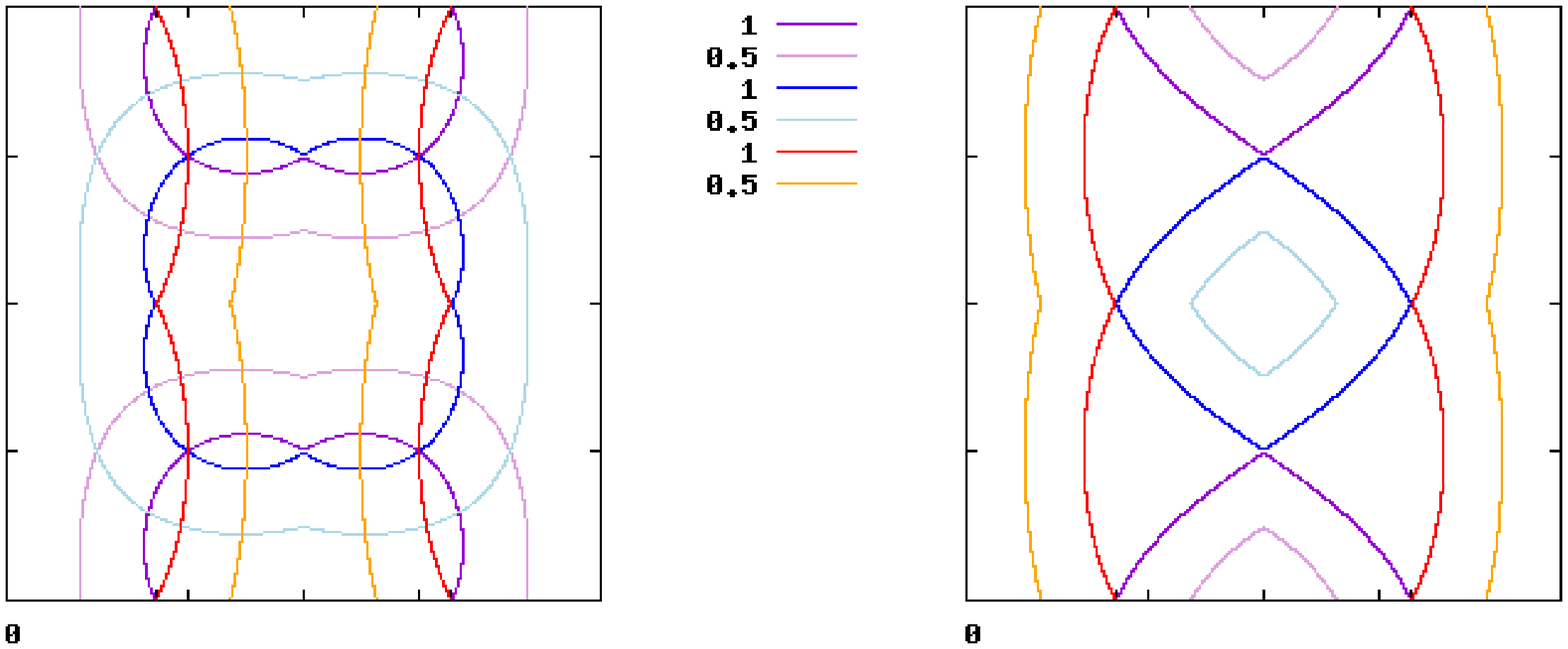}}
 \\
{\sc Figure 2A\hspace{6.3cm}Figure 2B}
\end{figure}
\end{center}
%g1(x,y)=sin(x)*sin(x)*cos(y)*cos(y)+abs(sin(x)*sin(x)*sin(y)*cos(y))+abs(sin(x)*cos(x)*cos(y))
%g2(x,y)=sin(x)*sin(x)*sin(y)*sin(y)+abs(sin(x)*sin(x)*sin(y)*cos(y))+abs(sin(x)*cos(x)*sin(y))
%g3(x,y)=cos(x)*cos(x)+abs(sin(x)*cos(x)*cos(y))+abs(sin(x)*cos(x)*sin(y))
%h1(x,y)=1-sin(x)*sin(x)*cos(y)*cos(y)+abs(sin(x)*sin(x)*sin(y)*cos(y))+abs(sin(x)*cos(x)*cos(y))
%h2(x,y)=1-sin(x)*sin(x)*sin(y)*sin(y)+abs(sin(x)*sin(x)*sin(y)*cos(y))+abs(sin(x)*cos(x)*sin(y))
%h3(x,y)=1-cos(x)*cos(x)+abs(sin(x)*cos(x)*cos(y))+abs(sin(x)*cos(x)*sin(y))
%set terminal png transparent crop size 640,480 xffffff \
%x000000 x202020 \
%xadd8e6 x9500d3 xdda0dd x0000ff xadd8e6 xff0000 xffa500
In Figure~2A we have plotted\footnote{ The figure has been obtained using Gnuplot 4.2 with a grid density of 1000 on both axes. For reference we have plotted the level lines for the value $0.5$, too. On the $\xi$-axis (horizontal) we have highlighted the values $\frac\pi4$, $\arctan\sqrt{2}$, $\frac{3\pi}4$ and $\pi-\arctan\sqrt{2}$. On the $\phi$-axis (vertical) we have highlighted the values $\frac{\pi}{4}$, $\frac\pi2$ and $\frac{3\pi}{4}$.} the level lines of the above functions  for the value $1$ (in violet, blue and red, respectively). This suggests that the ten parameter choices\\
\begin{minipage}{7cm}
\begin{itemize}
\item $(\xi,\phi)=(\arctan{\sqrt{2}},\frac{\pi}{4})$,
\item $(\xi,\phi)=(\pi-\arctan{\sqrt{2}},\frac{\pi}{4})$,
\item $(\xi,\phi)=(\arctan{\sqrt{2}},\frac{3\pi}{4})$,
\item $(\xi,\phi)=(\pi-\arctan{\sqrt{2}},\frac{3\pi}{4})$.
\end{itemize} 
\end{minipage}
\begin{minipage}{7cm}
\begin{itemize}
\item $(\xi,\phi)=(\frac{\pi}{4},\frac{\pi}{2})$,
\item $(\xi,\phi)=(\frac{3\pi}{4},\frac{\pi}{2})$,
\item $(\xi,\phi)=(\frac{\pi}{2},\frac{\pi}{4})$,
\item $(\xi,\phi)=(\frac{\pi}{2},\frac{3\pi}{4})$,
\item $(\xi,\phi)=(\frac{\pi}{4},0)$,
\item $(\xi,\phi)=(\frac{3\pi}{4},0)$,
\end{itemize} 
\end{minipage}\\
lead to an $L^\infty$-contractive semigroup -- as in fact can be checked directly.

Observe that $Y_{\xi,\phi}$ identifies Kirchhoff boundary conditions if and only if 
$$P_{Y_{\xi,\phi}}=\frac{1}{3}\begin{pmatrix}1 & 1 & 1\\ 1 & 1 & 1\\ 1 & 1 & 1\end{pmatrix},$$
i.e., if and only if $\xi=\arctan \sqrt{2}$ and $\phi=\frac{\pi}{4}$. A direct computation shows that the remaining nine cases correspond to boundary conditions defined by means of spaces $Y$ given by
$$\{(c,-c,-c):c\in\mathbb C\},\quad \{(c,c,-c):c\in\mathbb C\},\quad \{(c,-c,c):c\in\mathbb C\},$$ 
$$\{(c,c,0):c\in\mathbb C\},\quad \{(c,-c,0):c\in\mathbb C\},\quad \{(0,c,c):c\in\mathbb C\},$$ 
$$\{(0,c,-c):c\in\mathbb C\},\quad\{(c,0,c):c\in\mathbb C\},\quad \{(c,0,-c):c\in\mathbb C\}.$$ 
% that are ranges of the orthogonal projection matrices
% $$\frac{1}{3}\begin{pmatrix}1 & -1 & -1\\ -1 & 1 & 1\\ -1 & 1 & 1\end{pmatrix},\quad\frac{1}{3}\begin{pmatrix}1 & 1 & -1\\ 1 & 1 & -1\\ -1 & -1 & 1\end{pmatrix},
% \quad\frac{1}{3}\begin{pmatrix}1 & -1 & 1\\ -1 & 1 & -1\\ 1 & -1 & 1\end{pmatrix}\quad\hbox{and}$$
% $$\frac{1}{2}\begin{pmatrix}1 & 1 & 0\\ -1 & 1 & 1\\ -1 & 1 & 1\end{pmatrix},\quad\frac{1}{3}\begin{pmatrix}1 & 1 & -1\\ 1 & 1 & -1\\ -1 & -1 & 1\end{pmatrix},
% \quad\frac{1}{3}\begin{pmatrix}1 & -1 & 1\\ -1 & 1 & -1\\ 1 & -1 & 1\end{pmatrix}\quad\hbox{and}$$
% for $(\xi,\phi)=(\arctan{\sqrt{2}},\frac{3\pi}{4})$, $(\xi,\phi)=(\pi-\arctan{\sqrt{2}},\frac{\pi}{4})$ and $(\xi,\phi)=(\pi-\arctan{\sqrt{2}},\frac{3\pi}{4})$, respectively. 
While the last six subspaces only describe some decoupling of any of the three edges, we cannot find any physical interpretation for the first three boundary conditions. One can see that analogous boundary conditions give rise to $L^\infty$-contractive semigroups also in higher dimensional spaces $H={\mathbb C}^N$ for any $N\in\mathbb N$.

It ought to be remarked that not all relevant values become evident through the above plot: one can see that decoupled boundary conditions arise with $Y_{0,\phi}$ and $Y_{\frac{\pi}{2},\phi}$ for all $\phi\in[0,\pi)$ as well as with $Y_{\xi,0}$ and $Y_{\xi,\frac\pi2}$ for all $\xi\in [0,\pi)$. Hence, using again the computations performed in the case of $H={\mathbb C}^2$, we see that  $Y_{\frac{\pi}{2},\phi}$ lead to $L^\infty$-contractivity for $\phi\in\{0,\frac\pi4,\frac\pi2,\frac{3\pi}{4}\}$, and so do $Y_{\xi,0}$ and $Y_{\xi,\frac\pi2}$ for $\xi\in\{0,\frac\pi4,\frac\pi2,\frac{3\pi}{4}\}$ as well as $Y_{0,\phi}$ for all $\phi\in [0,\pi)$. We do not know whether further pairs $(\xi,\phi)$ leading to $L^\infty$-contractivity exist.

Moreover, a straightforward computation shows that $P_{Y_{\xi,\phi}}$ is a positive matrix if and only if $(\xi,\phi)\in [0,\frac{\pi}{2}]\times [0, \frac{\pi}{2}]$. Again, Kirchhoff boundary conditions are a singularity in a ``sea'' of non-submarkovian behaviours.

A similar procedure identifies all the $2$-dimensional subspaces of ${\mathbb C}^3$, i.e., all ranges of the matrices ${\rm Id}-P_{Y_{\xi,\phi}}$, $\xi,\phi\in [0,2\pi)$. However, plotting the level lines of the corresponding three functions (as we have done in Figure 2B in violet, blue and red, respectively), does not suggest any new pairs $(\xi,\phi)$ that lead to $L^\infty$-contractivity.

The general case of a semi-infinite star with arbitrarily (finitely) many edges can be treated likewise, using known formulae for hyperspherical coordinates.

 As already remarked, the above results carry over to case of dynamic boundary conditions and should be compared with the known properties of the heat equation with Wentzell--Robin boundary conditions in the scalar case, cf.~\cite{Bel94,MugRom07} and references therein. 

\section{Dynamic boundary conditions on the normal derivative}\label{dynneum}

In this section we consider a different setting by discussing a new kind of dynamics on the boundary. While the dynamic boundary conditions introduced in $\rm(AV)$ involve the trace, dynamic boundary conditions \emph{on the normal derivative} have also been considered in the literature, although less commonly (see~\cite{CasEngNag03,CarMug08}). Accordingly, the similar but different abstract initial-boundary value problem
\begin{equation}\tag{AVN}
\left\{\begin{array}{rcll}
\frac{\partial}{\partial t}{u}(t)&=& \Delta u(t), &t\geq 0,\\
\frac{\partial u}{\partial\nu } &\in& {\mathcal Y}, &t\geq 0,\\
\frac{\partial^2 }{\partial t\partial\nu }u(t)_{|\Gamma}&=&\delta P_{\mathcal Y}u_{|\Gamma}(t)+
P_{\mathcal Y}\left(\gamma \Delta_\Gamma -{\mathcal S}\right) \frac{\partial u}{\partial\nu }(t), &t\geq 0,\\
u(0)&=&u_0,\\
\frac{\partial u}{\partial \nu}(0)_{|\Gamma}&=&v_0,
\end{array}
\right.
\end{equation}
can be studied for $\gamma\in {\mathbb R}_+$ and ${\mathcal S}\in{\mathcal L}(L^2(\Omega;H))$. The parameter $\delta\in\mathbb C$ will be shown to influence the behaviour of the solutions to $\rm(AVN)$ in a curious way.

Consider a sesquilinear form $b_{\mathcal S}$ defined by
\begin{eqnarray*}
b_{\mathcal S} \left(\begin{pmatrix}f_1\\f_2\end{pmatrix},\begin{pmatrix}g_1\\ g_2\end{pmatrix}\right)&:=&\int_\Omega \left(\nabla f_1(x) | \nabla g_1(x)\right)_{H^n} dx-\delta\left(f_{1|\Gamma} | g_2\right)_{\mathcal Y} d\sigma(z)\\
&&\quad-  \left(f_2 | g_{1|\Gamma}\right)_{\mathcal Y} d\sigma(z) %+ \int_{\Gamma} ({\mathcal S}f_{|\Gamma}_1(z)|g_2(z))_{H} d\sigma(z)\\
 +\gamma\int_{\Gamma} (\nabla f_2(z)|\nabla g_2(z))_{H^{n-1}}d\sigma(z)+({\mathcal S}f_2|g_2)_{\mathcal Y},
\end{eqnarray*}
with dense domain 
$$W_{{\mathcal Y}}:=H ^1(\Omega;H)\times \left(H^s(\Gamma;H)\cap{\mathcal Y}\right),$$ 
where $s=0$ if $\gamma =0$ or $s=1$ if $\gamma>0$. Mimicking the proof of Theorem~\ref{formdir} we deduce a corresponding generation result (cf. also the discussion in~\cite[{\S}4.3]{CarMug08}).

\begin{theo}\label{formneu}
For any $\gamma\in {\mathbb R}_+$, $\delta\in\mathbb C$ and ${\mathcal S}\in{\mathcal L}(L^2(\Gamma;H))$ the sesquilinear form $b_{\mathcal Y}$ is continuous and elliptic (with respect to ${\mathcal L}^2$). The operator $B_{{\mathcal Y},{\mathcal S}}$ associated with $(b_{\mathcal S},W_{\mathcal Y})$ generates an analytic semigroup $(e^{tB_{\mathcal Y}})_{t\ge 0}$ on ${\mathcal L}^2$ with angle $\frac\pi2$. The semigroup is compact if and only if $\Omega,\Gamma$ have finite measure, provided that $H$ is finite dimensional.
Moreover, $b_{\mathcal Y}$ is accretive if $\delta=-1$ and ${\mathcal S}$ is accretive; it is symmetric if and only if $\delta=1$ and $\mathcal S$ is self-adjoint. In these cases the semigroups is contractive and self-adjoint, respectively.
\end{theo}

% \begin{proof}
% Continuity and ${\mathcal L}^2$-ellipticity can be checked as usual by the Cauchy--Schwarz inequality and simple perturbation arguments. Due to finite measure of $\Omega,\Gamma$, we can restrict ourselves to considering constant functions $f_1,f_2$ and deduce that the form $b_{\mathcal Y}$ is accretive if and only if
% $$-\delta{\rm Re}\int_{\Gamma} \left(f_1(z) | f_2(z)\right)_H d\sigma(z)-{\rm Re}\int_{\Gamma} \left(f_2(z) | f_1(z)\right)_H d\sigma(z)+ {\rm Re}({\mathcal S}f_2|f_2)_{\mathcal Y}\ge 0.$$
% The assertion follows by linearity.
% \end{proof}

With a proof similar to that of Proposition~\ref{identdir} we can show the following, see also~\cite[{\S}1.8]{Car08}.

\begin{prop}\label{identdir2}
Assume $\Omega$ to have $C^2$-boundary. For any $\gamma\in {\mathbb R}_+$, $\delta\in\mathbb C$ and ${\mathcal S}\in{\mathcal L}(L^2(\Omega;H))$ the operator $B_{{\mathcal Y},{\mathcal S}}$ associated with $(b_{\mathcal S},W_{\mathcal Y})$ is given by
\begin{eqnarray*}
D(B_{{\mathcal Y},{\mathcal S}})&=&\left\{\begin{pmatrix} f\\ g\end{pmatrix}\in W_{\mathcal Y}: \Delta u\in L^2(\Omega;H),\; \Delta u_{\Gamma}\in L^2(\Gamma;H),\hbox{ and }\frac{\partial f}{\partial \nu}\in L^2(\Gamma;H)\right\},\\
B_{{\mathcal Y},{\mathcal S}}&=&\begin{pmatrix} \Delta & 0\\
\delta P_{\mathcal Y}T & P_{\mathcal Y}\left(\gamma \Delta_\Gamma-{\mathcal S}\right)\end{pmatrix},
\end{eqnarray*}
where $T$ denotes the trace operator from $H^1(\Omega;H)$ to $H^\frac{1}{2}(\Gamma;H)$, cf.\ \cite[{\S}7.1]{CarMug09}.
\end{prop}

Thus, the semigroup associated with $B_{{\mathcal Y},{\mathcal S}}$ yields the solution to the abstract initial-boundary value problem
\begin{equation*}
\left\{\begin{array}{rcll}
\frac{\partial{u}}{\partial t}(t)&=& \Delta u(t), &t\geq 0,\\
\frac{\partial u}{\partial\nu}(t)&\in &\mathcal Y,&t\ge 0,\\
\frac{\partial^2 u}{\partial t\partial \nu}(t)&=&\delta \left(P_{\mathcal Y}u(t)\right)+\left(\left(\gamma \Delta_\Gamma -{\mathcal S}\right)\frac{\partial u}{\partial \nu}(t)\right), &t\ge 0,\\
u(0)&=&u_{0},\\
\frac{\partial u}{\partial\nu}(0)&=& w_{0}.
\end{array}
\right.
\end{equation*}

Ouhabaz's criterion may be promptly applied to this setting, too. We omit the easy proof.

\begin{prop}
Impose Assumptions~\ref{locality}. Let $C_H$ be a closed subspace or a closed order interval of $H$. Consider the closed convex subsets $C_{L^2(\Omega;H)}$ and $C_{\mathcal Y}$. Then $(e^{tB_{{\mathcal Y},{\mathcal S}}})$ leaves invariant $C_{L^2(\Omega;H)}\times C_{\mathcal Y}$ if and only if the compatibility condition
\begin{eqnarray*}
\delta{\rm Re}\left(P_{C_{L^2(\Omega;H)}}f_{1|\Gamma} | \left(I-P_{C_{\mathcal Y}}\right) f_2\right)_{\mathcal Y}+{\rm Re} \left(P_{C_{\mathcal Y}} f_2 | \left(I-P_{C_{L^2(\Omega;H)}}\right) f_{1|\Gamma}\right)_{\mathcal Y} %+ \int_{\Gamma} ({\mathcal S}f_{|\Gamma}_1(z)|g_2(z))_{H} d\sigma(z)\\
\le{\rm Re}\left(SP_{C_{\mathcal Y}}f_2|\left(I-P_{C_{\mathcal Y}}\right)f_2\right)_{\mathcal Y}
\end{eqnarray*}
holds for all $f_1\in H^1(\Omega;H)$ and all $f_2\in {\mathcal Y}$.
\end{prop}

\begin{exa}
Impose Assumptions~\ref{locality}. Then, by linearity $(e^{tB_{Y,S}})$ is positive if and only if $\delta=1$ and 
$${\rm Re}\left(SP_{D_{\mathcal Y}}f_2|\left(I-P_{D_{\mathcal Y}}\right)f_2\right)_{\mathcal Y}\ge 0,$$
i.e., if and only if $\delta=1$ and the semigroup on $H$ generated by $-S$ is positive.
\end{exa}

\begin{rem}
It is easy to see that by similar methods one can also treat the parabolic problem
\begin{equation*}%\tag{VVD$_3$}
\left\{\begin{array}{rcll}
\frac{\partial}{\partial t}{u}(t)&=& \Delta u(t), &t\geq 0,\\
\frac{\partial}{\partial t} P_{\mathcal Y} u(t)_{|\partial \Omega}&=&-P_{\mathcal Y}\frac{\partial u(t)}{\partial \nu}+R_1 u(t)_{|\partial \Omega}, &t\geq 0,\\
\frac{\partial}{\partial t} P_{{\mathcal Y}^\perp}\frac{\partial u(t)}{\partial \nu}&=&-P_{{\mathcal Y}^\perp}u(t)_{|\partial \Omega}+R_2\frac{\partial u(t)}{\partial \nu}, &t\geq 0,\\
u(0)&=&u_{0},&\\
u(0)_{|\Gamma}&=&u_{1},&\\
\frac{\partial_D u}{\partial \nu}(0)_{|\Gamma}&=& u_2,&\\
\end{array}\right.
\end{equation*}
for some $R_1\in{\mathcal L}(H^1(\Omega;H),\mathcal Y)$ and $R_2\in{\mathcal L}(H^1(\Omega;H),\mathcal Y^\perp)$.  In this case the state space is $L^2(\Omega;H)\times {\mathcal Y}\times {\mathcal Y}^\perp$. We omit the details.
\end{rem}


\begin{thebibliography}{10}

\bibitem{AgmDouNir64}
S.~Agmon, A.~Douglis, and L.~Nirenberg.
\newblock Estimates near the boundary for solutions of elliptic partial
  differential equations satisfying general boundary conditions {II}.
\newblock {\em Comm. Pure Appl. Math.}, 17:35--92, 1964.

\bibitem{AlbCacFic07}
S.~Albeverio, C.~Cacciapuoti, and D.~Finco.
\newblock Coupling in the singular limit of thin quantum waveguides.
\newblock {\em J. Math. Phys.}, 48:032103, 2007.

\bibitem{AliNic93}
F.~{Ali Mehmeti} and S.~Nicaise.
\newblock Nonlinear interaction problems.
\newblock {\em Nonlinear Anal., Theory Methods Appl.}, 20:27--61, 1993.

\bibitem{Ama01}
H.~Amann.
\newblock Elliptic operators with infinite-dimensional state spaces.
\newblock {\em J. Evol. Equ.}, 1:143--188, 2001.

\bibitem{AliBelNic01}
F.~{Ali Mehmeti} and{J. von} Below and S.~Nicaise, editors.
\newblock {\em Partial Differential Equations on Multistructures (Proc. Luminy
  1999)}, volume 219 of {\em Lect. Notes Pure Appl. Math.}
\newblock Marcel Dekker, New York, 2001.

\bibitem{Are06}
W.~Arendt.
\newblock Heat {K}ernels -- {M}anuscript of the $9^{\rm th}$ {I}nternet
  {S}eminar, 2006.
\newblock Freely available at
  {\url{http://www.uni-ulm.de/fileadmin/website_uni_ulm/mawi.inst.020/arendt/d%
ownloads/internetseminar.pdf}}.

\bibitem{AreBatHie01}
W.~Arendt, C.J.K. Batty, M.~Hieber, and F.~Neubrander.
\newblock {\em Vector-{V}alued {L}aplace {T}ransforms and {C}auchy {P}roblems},
  volume~96 of {\em Monographs in Mathematics}.
\newblock Birkh{\"a}user, Basel, 2001.

\bibitem{AreBieEls08}
W.~Arendt, M.~Biegert, and T.~{ter Elst}.
\newblock Diffusion determines the manifold.
\newblock \url{http://arxiv.org/abs/0806.0437}.

\bibitem{AreMetPal03}
W.~Arendt, G.~Metafune, D.~Pallara, and S.~Romanelli.
\newblock The {L}aplacian with {W}entzell--{R}obin boundary conditions on
  spaces of continuous functions.
\newblock {\em Semigroup Forum}, 67:247--261, 2003.

\bibitem{Bel94}
{J. von} Below.
\newblock {\em {Parabolic Networks Equations}}.
\newblock T{\"u}binger Universit{\"a}tsverlag, T{\"u}bingen, 1994.

\bibitem{BelMug10}
{J. von} Below and D.~Mugnolo.
\newblock Spectral asymptotics for diffusive interface problems.
\newblock Preprint, 2010.

\bibitem{BelNic96}
{J. von} Below and S.~Nicaise.
\newblock Dynamical interface transition in ramified media with diffusion.
\newblock {\em Comm. Partial Differ. Equations}, 21:255--279, 1996.

\bibitem{Bre83}
H.~Brezis.
\newblock {\em Analyse {F}onctionelle - {T}h{\'e}orie et {A}pplications}.
\newblock Masson, Paris, 1983.

\bibitem{Cam80}
H.~Camerer.
\newblock {\em Die Elektrotonische Spannungsbreitung im Soma, Dendritenbaum und
  Axon von Nervenzellen}.
\newblock PhD thesis, Universit{\"a}t T{\"u}bingen, 1980.

\bibitem{Car08}
S.~Cardanobile.
\newblock {\em Diffusion Systems and Heat Equations on Networks}.
\newblock PhD thesis, Universit{\"a}t Ulm, 2008.

\bibitem{CarMug08}
S.~Cardanobile and D.~Mugnolo.
\newblock Qualitative properties of coupled parabolic systems of evolution
  equations.
\newblock {\em Ann. Sc. Norm. Super. Pisa, Cl. Sci., V Ser.}, 2:287--312, 2008.

\bibitem{CarMug09}
S.~Cardanobile and D.~Mugnolo.
\newblock Parabolic systems with coupled boundary conditions.
\newblock {\em J. Differ. Equ.}, 247:1229--1248, 2009.

\bibitem{CarMug09b}
S.~Cardanobile and D.~Mugnolo.
\newblock Towards a gauge theory for evolution equations on vector-valued
  spaces.
\newblock {\em J. Math. Phys.}, (in press), 2009.

\bibitem{CarMugNit08}
S.~Cardanobile, D.~Mugnolo, and R.~Nittka.
\newblock Well-posedness and symmetries of strongly coupled network equations.
\newblock {\em J. Phys. A}, 41:055102, 2008.

\bibitem{CasEngNag03}
V.~Casarino, K.-J. Engel, R.~Nagel, and G.~Nickel.
\newblock A semigroup approach to boundary feedback systems.
\newblock {\em Int. Equations Oper. Theory}, 47:289--306, 2003.

\bibitem{CocFavGol08}
G.M. Coclite, A.~Favini, G.R. Goldstein, J.A. Goldstein, and S.~Romanelli.
\newblock Continuous dependence on the boundary conditions for the {W}entzell
  {Laplacian}.
\newblock {\em Sem. Forum}, 1:101--108, 2008.

\bibitem{Cro07}
M.~Crouzeix.
\newblock Numerical range and functional calculus in {H}ilbert space.
\newblock {\em J. Funct. Anal.}, 244:668--690, 2007.

\bibitem{DauLio88}
R.~Dautray and J.-L. Lions.
\newblock {\em {Mathematical Analysis and Numerical Methods for Science and
  Technology, Vol. 2}}.
\newblock Springer-Verlag, Berlin, 1988.

\bibitem{DenHiePru03}
R.~Denk, M.~Hieber, and J.~Pr{\"u}ss.
\newblock {\em ${{R}}$-{B}oundedness, {F}ourier {M}ultipliers and {P}roblems of
  {E}lliptic and {P}arabolic {T}ype}, volume 788 of {\em Mem. Am. Math. Soc.}
\newblock Amer. Math. Soc., Providence, RI, 2003.

\bibitem{DesSch84}
W.~Desch and W.~Schappacher.
\newblock On relatively bounded perturbations of linear ${C}_0$-semigroups.
\newblock {\em Ann. Sc. Norm. Super. Pisa, Cl. Sci.}, 11:327--341, 1984.

\bibitem{ExnPos05}
P.~Exner and O.~Post.
\newblock Convergence of spectra of graph-like thin manifolds.
\newblock {\em J. Geom. Phys.}, 54:77--115, 2005.

\bibitem{FavGolGol02}
A.~Favini, G.R. Goldstein, J.A. Goldstein, and S.~Romanelli.
\newblock The heat equation with generalized {W}entzell boundary condition.
\newblock {\em J. Evol. Equ.}, 2:1--19, 2002.

\bibitem{FulKucWil07}
S.~A. Fulling, P.~Kuchment, and J.~H. Wilson.
\newblock Index theorems for quantum graphs.
\newblock {\em J. Phys. A}, 40:14165--14180, 2007.

\bibitem{GilTru01}
D.~Gilbarg and N.~Trudinger.
\newblock {\em {Elliptic Partial Differential Equations of Second Order}}.
\newblock Classics in Mathematics. Springer-Verlag, Berlin, 2001.

\bibitem{Gra08}
L.~Grafakos.
\newblock {\em {Classical Fourier Analysis}}, volume 249 of {\em Graduate Texts
  in Mathematics}.
\newblock Springer-Verlag, Berlin, 2008.

\bibitem{HanZua95}
S.~Hansen and E.~Zuazua.
\newblock {Exact controllability and stabilization of a vibrating string with
  an interior point mass}.
\newblock {\em SIAM J. Control Optimization}, 33:1357--1391, 1995.

\bibitem{ItoKap02}
K.~Ito and F.~Kappel.
\newblock {\em {Evolution Equations and Approximations}}, volume~61 of {\em
  Adv. Math. Appl. Sci.}
\newblock World Scientific, 2002.

\bibitem{JohLin01}
W.B. Johnson and J.~{J. Lindenstrauss}, editors.
\newblock {\em {Handbook of the geometry of {B}anach spaces. Vol. 1}}.
\newblock Elsevier, Amsterdam, 2001.

\bibitem{KanKlaVoi09}
U.~Kant, T.~Klau{\ss}, J.~Voigt, and M.~Weber.
\newblock Dirichlet forms for singular one-dimensional operators and on graphs.
\newblock {\em J. Evol. Equ.}, 9:(in press), 2009.

\bibitem{Kat66}
T.~Kato.
\newblock {\em {Perturbation Theory for Linear Operators}}.
\newblock Classics in Mathematics. Springer-Verlag, New York, 1966.

\bibitem{KosPotSch08}
V.~Kostrykin, J.~Potthoff, and R.~Schrader.
\newblock Contraction semigroups on metric graphs.
\newblock In P.~Exner, J.~Keating, P.~Kuchment, T.~Sunada, and A.~Teplyaev,
  editors, {\em {Analysis on Graphs and its Applications}}, volume~77 of {\em
  Proceedings of Symposia in Pure Mathematics}, pages 423--458, Providence, RI,
  2008. Amer. Math. Soc.

\bibitem{KosSch99}
V.~Kostrykin and R.~Schrader.
\newblock Kirchhoff's rule for quantum wires.
\newblock {\em J. Phys. A}, 32:595--630, 1999.

\bibitem{Kuc04}
P.~Kuchment.
\newblock Quantum graphs {I}: {S}ome basic structures.
\newblock {\em Waves in Random Media}, 14:107--128, 2004.

\bibitem{Kuc08}
P.~Kuchment.
\newblock Quantum graphs: an introduction and a brief survey.
\newblock In P.~Exner, J.~Keating, P.~Kuchment, T.~Sunada, and A.~Teplyaev,
  editors, {\em {Analysis on Graphs and its Applications}}, volume~77 of {\em
  Proceedings of Symposia in Pure Mathematics}, pages 291--314, Providence, RI,
  2008. Amer. Math. Soc.

\bibitem{KucZen03}
P.~Kuchment and H.~Zeng.
\newblock Asymptotics of spectra of {N}eumann {L}aplacians in thin domains.
\newblock In Y.~Karpeshina {et al.}, editor, {\em {Advances in Differential
  Equations and Mathematical Physics (Proc. Birmingham 2002)}}, volume 327 of
  {\em Contemp. Math.}, pages 199--213, Providence, 2003. Amer. Math. Soc.

\bibitem{LagLeuSch94}
J.E. Lagnese, G.~Leugering, and E.J.P.G. Schmidt.
\newblock {\em Modeling, {A}nalysis, and {C}ontrol of {D}ynamic {E}lastic
  {M}ulti-{L}ink {S}tructures}.
\newblock Systems and Control: Foundations and Applications. Birkh{\"a}user,
  Basel, 1994.

\bibitem{LioMag72}
J.L. Lions and E.~Magenes.
\newblock {\em Non-{H}omogeneous {B}oundary {V}alue {P}roblems and
  {A}pplications}, volume 181--183 of {\em Grundlehren der mathematischen
  Wissenschaften}.
\newblock Springer-Verlag, Berlin, 1972.

\bibitem{Lum80a}
G.~Lumer.
\newblock Espaces ramifi{\'e}s et diffusion sur les r{\'e}seaux topologiques.
\newblock {\em C.R. Acad. Sc. Paris}, 291:627--630, 1980.

\bibitem{MajEvaJac93}
G.~Major, J.D. Evans, and J.J. Jack.
\newblock Solutions for transients in arbitrarily branching cables: I. voltage
  recording with a somatic shunt.
\newblock {\em Biophys. J.}, 65:423--449, 1993.

\bibitem{ManVogVoi05}
A.~Manavi, H.~Vogt, and J.~Voigt.
\newblock {Domination of semigroups associated with sectorial forms}.
\newblock {\em Journal of Operator Theory}, 54:9--25, 2005.

\bibitem{Maz85}
V.G. Maz'ya.
\newblock {\em Sobolev {S}paces}.
\newblock Springer-Verlag, Berlin, 1985.

\bibitem{Mey91}
P.~Meyer-Nieberg.
\newblock {\em {B}anach Lattices}.
\newblock Universitext. Springer-Verlag, Berlin, 1991.

\bibitem{MorFisPes08}
Y.~Mori, G.I. Fishman, and C.S. Peskin.
\newblock Ephaptic conduction in a cardiac strand model with {3D}
  electrodiffusion.
\newblock {\em PNAS}, 105:6463--6468, 2008.

\bibitem{Mug07}
D.~Mugnolo.
\newblock Gaussian estimates for a heat equation on a network.
\newblock {\em Networks Het. Media}, 2:55--79, 2007.

\bibitem{Mug08c}
D.~Mugnolo.
\newblock Asymptotics of semigroups generated by operator matrices.
\newblock \url{http://arxiv.org/abs/0801.1963}, 2008.

\bibitem{Mug08}
D.~Mugnolo.
\newblock A variational approach to strongly damped wave equations.
\newblock In H.~Amann {et al.}, editor, {\em Functional Analysis and Evolution
  Equations -- The G{\"u}nter Lumer Volume}, pages 503--514. Birkh{\"a}user,
  Basel, 2008.

\bibitem{MugNit09}
D.~Mugnolo and R.~Nittka.
\newblock Properties of representations of operators acting between spaces of
  vector-valued functions.
\newblock \url{http://arxiv.org/abs/0903.2038}, 2009.

\bibitem{MugNitPos10}
D.~Mugnolo, R.~Nittka and O.~Post.
\newblock Convergence of sectorial operators on varying {H}ilbert spaces.
\newblock \url{(Submitted)}

\bibitem{MugRom06}
D.~Mugnolo and S.~Romanelli.
\newblock Dirichlet forms for general {W}entzell boundary conditions, analytic
  semigroups, and cosine operator functions.
\newblock {\em Electronic J. Differ. Equ.}, 118:1--20, 2006.

\bibitem{MugRom07}
D.~Mugnolo and S.~Romanelli.
\newblock Dynamic and generalized {W}entzell node conditions for network
  equations.
\newblock {\em Math. Meth. Appl. Sci.}, 30:681--706, 2007.

\bibitem{Nag86}
R.~Nagel, editor.
\newblock {\em One-{P}arameter {S}emigroups of {P}ositive {O}perators}, volume
  1184 of {\em Lect. Notes Math.}
\newblock Springer-Verlag, Berlin, 1986.

\bibitem{Nic85}
S.~Nicaise.
\newblock Some results on spectral theory over networks, applied to nerve
  impulse transmission.
\newblock In C.~Brezinsky {et al.}, editor, {\em Polyn\^{o}mes Orthogonaux et
  Applications (Proc. Bar-le-Duc 1984)}, volume 1171 of {\em Lect. Notes.
  Math.}, pages 532--541, Berlin, 1985. Springer-Verlag.

\bibitem{Nic88}
S.~Nicaise.
\newblock Elliptic operators on elementary ramified spaces.
\newblock {\em Int. Equations Oper. Theory}, 11:230--257, 1988.

\bibitem{Ouh05}
E.M. Ouhabaz.
\newblock {\em Analysis of {H}eat {E}quations on {D}omains}, volume~30 of {\em
  LMS Monograph Series}.
\newblock Princeton University Press, Princeton, 2005.

\bibitem{Pos06}
O.~Post.
\newblock Spectral convergence of quasi-one-dimensional spaces.
\newblock {\em Ann. Henri Poincare}, 7:933--973, 2006.

\bibitem{Pos07}
O.~Post.
\newblock First order operators and boundary triples.
\newblock {\em Russ. J. Math. Phys.}, 14:482--492, 2007.

\bibitem{Ral59}
W.~Rall.
\newblock Branching dendritic trees and motoneurone membrane resistivity.
\newblock {\em Exp. Neurol.}, 1:491--527, 1959.

\bibitem{Sho94}
R.E. Showalter.
\newblock {\em {H}ilbert Space Methods for Partial Differential Equations}.
\newblock Electronic Journal of Differential Equations, San Marcos, TX, 1994.

\bibitem{Sho97}
R.E. Showalter.
\newblock {\em Monotone Operator in {B}anach Space and Partial Differential
  Equations}, volume~49 of {\em Math. Surveys and Monographs}.
\newblock Amer. Math. Soc., Providence, RI, 1997.

\bibitem{VazVit08}
J.~L. V{\'a}zquez and E.~Vitillaro.
\newblock Heat equation with dynamical boundary conditions of reactive type.
\newblock {\em Comm. Partial Differ. Equations}, 33:561--612, 2008.

\end{thebibliography}
\end{document}